\documentclass[12pt]{amsart}
\usepackage{amssymb,amsthm,amsmath,amscd,rsfs}
\usepackage{latexsym}

\newcommand{\NN}{{\mathbb{N}}}
\newcommand{\ZZ}{{\mathbb{Z}}}
\newcommand{\QQ}{{\mathbb{Q}}}
\newcommand{\RR}{{\mathbb{R}}}
\newcommand{\CC}{{\mathbb{C}}}

\newcommand{\KK}{{\mathbb{K}}}

\newcommand{\XXX}{{\mathscr{X}}}
\newcommand{\YYY}{{\mathscr{Y}}}

\newcommand{\UUU}{{\mathscr{U}}}

\newcommand{\BBB}{{\mathscr{B}}}
\newcommand{\EEE}{{\mathscr{E}}}
\newcommand{\AAA}{{\mathscr{A}}}
\newcommand{\TTT}{{\mathscr{T}}}
\newcommand{\CCC}{{\mathscr{C}}}

\newcommand{\MMM}{{\mathscr{M}}}
\newcommand{\LLL}{{\mathscr{L}}}

\newcommand{\OO}{{\mathcal{O}}}

\newcommand{\codim}{\operatorname{codim}}
\newcommand{\ch}{\operatorname{char}}

\newcommand{\Ker}{\operatorname{Ker}}

\newcommand{\Spec}{\operatorname{Spec}}
\newcommand{\Spf}{\operatorname{Spf}}

\newcommand{\Alb}{\operatorname{Alb}}

\newcommand{\cherncl}{{c}}

\newcommand{\ord}{\operatorname{ord}}

\newcommand{\Tr}{\operatorname{Tr}}

\newcommand{\aff}{\mathrm{aff}}

\newcommand{\an}{\mathrm{an}}
\newcommand{\trop}{\mathrm{trop}}
\newcommand{\val}{\operatorname{val}}
\newcommand{\Val}{\operatorname{Val}}
\newcommand{\str}{\operatorname{str}}
\newcommand{\red}{\operatorname{red}}

\theoremstyle{plain}
\newtheorem{Theorem}{Theorem}[section]
\newtheorem{Lemma}[Theorem]{Lemma}
\newtheorem{Proposition}[Theorem]{Proposition}
\newtheorem{Corollary}[Theorem]{Corollary}
\newtheorem{Main-Theorem}[Theorem]{Main Theorem}
\newtheorem{Theorem-Definition}[Theorem]{Theorem-Definition}
\theoremstyle{definition}
\newtheorem{Definition}[Theorem]{Definition}
\newtheorem{Remark}[Theorem]{Remark}
\newtheorem{Conjecture}[Theorem]{Conjecture}
\newtheorem{Claim}{Claim}

\newtheorem{Key Fact}[Theorem]{Key Fact}
       
\newtheorem*{n-c}{Notation and convention}      
\newtheorem{citeTheorem}[Theorem]{Theorem}

\renewcommand{\theTheorem}{\arabic{section}.\arabic{Theorem}}
\renewcommand{\theClaim}{\arabic{section}.\arabic{Theorem}.\arabic{Claim}}
\renewcommand{\theequation}{\arabic{section}.\arabic{Theorem}.\arabic{Claim}}

\newtheorem{art}[Theorem]{}

\newcommand{\Acal}{{\mathscr A}}

\newcommand{\Bcal}{{\mathscr B}}

\newcommand{\Ccal}{{\mathscr C}}

\newcommand{\Dcal}{{\mathscr D}}

\newcommand{\Ecal}{{\mathscr E}}

\newcommand{\Xcal}{{\mathscr X}}

\newcommand{\Ycal}{{\mathscr Y}}

\newcommand{\cdop}{{\mathbb C}}

\newcommand{\rdop}{{\mathbb R}}

\newcommand{\kdop}{{\mathbb K}}

\newcommand{\rtor}{{\rdop^n/\Lambda}}
\newcommand{\Deltabar}{{\overline{\Delta}}}
\newcommand{\Ccalbar}{{\overline{\Ccal}}}

\newcommand{\Sigmabar}{{\overline{\Sigma}}}

\newcommand{\valbar}{{\overline{\val}}}
\newcommand{\Xan}{{X^{\rm an}}}

\newcommand{\Aan}{{A^{\rm an}}}
\newcommand{\relint}{{\rm relint}}
\newcommand{\kcirc}{{\rm \kdop^\circ}}

\begin{document}

\title
[Geometric Bogomolov conjecture]
{Geometric Bogomolov conjecture for abelian varieties
\\
and some results for those 
with some degeneration
(with an appendix by Walter Gubler: The minimal dimension  
of a canonical measure)}
\author
{Kazuhiko Yamaki}
\email{yamaki.kazuhiko.6r@kyoto-u.ac.jp}
\date{
October 9, 2012 (Version 7.3).}
\address
{Institute for the Promotion of Excellence in Higher Education,
Kyoto University, Kyoto, 606-8501, Japan}

\begin{abstract}
In this paper, we formulate the geometric Bogomolov conjecture for abelian varieties, and give some partial answers to it.
In fact, we insist in a main theorem that
under some degeneracy condition,
a closed subvariety of an abelian variety does not have a dense subset of
small points
if
it is a non-special subvariety.
The key of the proof is the study of the minimal dimension of the components of a canonical measure on the tropicalization of the closed subvariety. Then we can apply the tropical version of equidistribution theory due to Gubler. This article includes an appendix by Walter Gubler. 
He shows that the minimal dimension of the components of a canonical  
measure is equal to the dimension of the abelian part of the subvariety.
We can apply this result to make a further contribution to the geometric
Bogomolov conjecture.
\end{abstract}

\maketitle


\section*{Introduction}

\subsection{Motivation and statements}

Let $K$ be a number field, 
or a function field 
over a base field $k$.
We fix an algebraic closure $\overline{K}$ of $K$.
Let $A$ be an  abelian variety over $\overline{K}$ and
let $L$ be an ample line bundle on $A$,
and assume it is even,
i.e., $[-1]^{\ast} L = L$.
Then the 
canonical height
function $\hat{h}_{L}$ associated with $L$,
also called the N\'eron-Tate height,
is a 
semi-positive definite quadratic form on 
$A 
\left(
\overline{K}
\right)$.
It is  well known that
$\hat{h}_{L} (x) = 0$ if $x$ is a torsion point.
Let $X$ be a closed subvariety of 
$A$.
We
put 
\[
X (\epsilon ; L)
:=
\left\{
x \in X 
\left(
\overline{K}
\right)
\left|
\hat{h}_{L} (x) \leq \epsilon
\right.
\right\}
\]
for
a positive real number $\epsilon > 0$.
Then the Bogomolov conjecture for abelian varieties insists that
there should exist an $\epsilon > 0$ such that
$X (\epsilon ; L)$
is not Zariski dense in $X$,
unless $X$ is a kind of ``exceptional'' closed subvarieties,
such as torsion subvarieties for example.

In the case where $K$ is a number field,
namely, in the arithmetic case,
this conjecture was solved more than ten years ago,
known as
a theorem of Zhang:

\begin{citeTheorem}[Corollary 3 in \cite{zhang2},
arithmetic version of
Bogomolov conjecture for abelian varieties]
\label{ABC}
Let $K$ be a number field.
If $X$ is not a torsion subvariety,
then there is an $\epsilon > 0$ such that
$X (\epsilon ; L)$
is not Zariski dense in $X$.
\end{citeTheorem}

The Bogomolov conjecture is originally
a statement concerning the
jacobian of a curve and an embedding of the curve,
that is,
$A$ is a jacobian and $X$ is an embedded curve.
It is called the Bogomolov conjecture
for curves,
and
is proved by Ullmo in \cite{ullmo}
in case that $K$ is a number field,
at the same time
when Zhang proved Theorem~\ref{ABC}.
The ideas of Ullmo and Zhang are same;
they are based 
on equidistribution theory.
The outline of this argument
will be recalled later in this introduction.

When $K$ is a finitely generated field over $\QQ$,
a kind of arithmetic height functions can be defined
after a choice of
polarizations of $K$, due to  Moriwaki \cite{moriwaki5}.
It is still an arithmetic setting namely,
and the Bogomolov conjecture for abelian varieties
with respect to the height
associated with a big polarization
has been proved
by Moriwaki himself.
The classical geometric height
is also a kind of Moriwaki's arithmetic height,
but
it does not arise from a big polarization --- rather a degenerate one.
Hence we cannot say anything about the geometric version
of the conjecture
with Moriwaki's theory.

How about the geometric case,
that is, the case where $K$ is a function field 
over an algebraically
closed field $k$
and
the height is the classical geometric height?
In this case,
we cannot expect the same statement as Theorem~\ref{ABC}
because a subvariety defined over the constant field
can have a dense subset of points of height zero.
Accordingly, we have to
reformulate the conjecture,
or have to consider it in a restricted situation.

The Bogomolov conjecture for curves
over a function field of one variable
has been studied for a long time
as one of the important special case 
of the conjecture over a function field.
In this case,
the ``exceptional'' ones are the isotrivial curves.
In characteristic $0$,
Cinkir proved this version of conjecture in \cite{cinkir}.
In positive characteristic,
the conjecture for curves is still open,
but there are some partial answers 
such as
in
\cite{moriwaki3} by Moriwaki and
in
\cite{yamaki1,yamaki2} by the author.

Another important result
on the Bogomolov conjecture in the geometric setting
is the one due to Gubler.
He proved in \cite{gubler2}
the following theorem:
\begin{citeTheorem} [Theorem 1.1 of \cite{gubler2}]
\label{thmofgubler}
Assume that there is a place $v$ at which the abelian variety
$A$ is totally degenerate.
Then $X (\epsilon ; L)$ is not Zariski dense in $X$ for some
$\epsilon > 0$
unless $X$ is a torsion subvariety.
\end{citeTheorem}
In this theorem,
the ``exceptional'' subvarieties are
the torsion subvarieties,
same
as in the arithmetic case,
because
there do not appear constant subvarieties
in this setting.

In this paper,
we discuss the Bogomolov conjecture in the geometric setting.
This paper has two goals: One is to give a precise formulation of
the geometric version of the Bogomolov conjecture 
for arbitrary abelian varieties. This may be well known to the
experts but seems to be lacking in the literature. 
The other is to prove
the conjecture 
under a certain condition on degeneration.
It is much more
general than the case of totally degenerate abelian varieties considered
in Theorem~\ref{thmofgubler}.

Let us give more details of each goal of ours.
Let $X$ be an irreducible closed subvarieties of $A$.
Let $\left( A^{\overline{K}/k} , \Tr_{A}^{\overline{K}/k} \right)$ denote the
$\overline{K}/k$-trace
of $A$ (cf. \S~\ref{appendix}).
Then $X$ is said to be
\emph{special} 
if there are a torsion point $\tau \in A ( \overline{K})$
and a closed subvariety $X'$ of $A^{\overline{K}/k}$ such that
\[
X = G_{X} + \Tr_{A}^{\overline{K}/k} \left( X'_{\overline{K}} \right) +\tau
,
\]
where $G_{X}$ is the stabilizer of $X$
(cf. \S~\ref{specialsub}).
Note that if there is a place $v$ at which $A$ is totally degenerate,
then the notion of special subvarieties coincides with
that of torsion subvarieties 
since the $\overline{K}/k$-trace is trivial.
We will see that
any special subvariety has dense small points
(cf. Corollary~\ref{special-small}).
Our geometric Bogomolov conjecture insists that 
the converse should hold true:

\begin{Conjecture} [cf. Conjecture~\ref{GBC}
]
Let $K$ be a function field.
Let all $A$, $L$ and $X$
be as above.
Then there exists an $\epsilon > 0$
such that $X (\epsilon ; L)$ is not Zariski dense
in $X$ unless $X$ is a special subvariety.
\end{Conjecture}

For an irreducible closed subvariety
$X \subset A$
and a place $v$ of $\overline{K}$,
we can define an integer $b (X_{v})$,
called the dimension of abelian part of $X_{v}$
(cf. \S~\ref{abpt}).
We do not give its definition here
because it is a little bit complicated,
but let us describe what it is
in the case where
$A_{v}$ is the product of an abelian variety
$B$ with good reduction and a totally degenerate
abelian variety;
if $\alpha : A_{v} \to B$ is the projection,
then $b(X_{v})$ coincides with $\dim \alpha (X)$.

We can see that if
there is a place $v$ with 
$\dim (X / G_{X}) > b ((X / G_{X})_{v})$,
then $X$ is not a special subvariety
(cf. Proposition~\ref{compprop} (1)).
Hence if our conjecture holds true,
then such an $X ( \epsilon ; L)$ should not be dense
in $X$ for some $\epsilon > 0$.
In fact,
we will show the following result
as our main theorem of this paper:

\begin{Theorem} [cf. Theorem~\ref{maintheorem}
]
\label{thoremint}
Assume that there exists a place $v$ such that
$\dim (X / G_{X}) > b ((X / G_{X})_{v})$.
Then $X (\epsilon ; L)$ is not Zariski dense 
for some $\epsilon > 0$.
\end{Theorem}

Note that
this theorem leads us to a generalization of Theorem~\ref{thmofgubler}.
In fact,
if $A_{v}$ is totally degenerate for some place $v$
and $X$ is a non-torsion subvariety,
then we see $b \left( A_{v} \right) = 0$ and hence
$b \left( (X / G_{X})_{v} \right) = 0$ by Lemma~\ref{quotientb}.
Therefore
$X ( \epsilon ; L)$ is not dense in $X$ for some $\epsilon > 0$
by Theorem~\ref{thoremint}
if $\dim \left( X / G_{X} \right) > 0$,
and
by Lemma~\ref{zerodimensional}
if $\dim \left( X / G_{X} \right) = 0$.

\subsection{Ideas}

We
would like here to describe the idea of our proof.
To do that,
let us recall the proof of Theorem~\ref{ABC}
and that of Theorem~\ref{thmofgubler},
which gives us a basic strategy.

First we recall the admissible metric.
Let $A$ be an abelian variety over $\CC$
and
let
$L$ be an even ample line bundle on $X$.
It is well known that there is a canonical
hermitian metric $h_{can}$ on $L$,
called the canonical metric,
such that
$[n]^{\ast} \cherncl_{1} (L , h_{can}) = 
n^{2} \cherncl_{1} (L , h_{can})$
and that
$\cherncl_{1} (L , h_{can})$ is smooth and positive,
where $\cherncl_{1} (L , h_{can})$ is the curvature form.
For a closed subvariety $X \subset A$
of dimension $d$,
put
\[
\mu_{X,L} := 
\frac{1}{\deg_{L} (X)}
\cherncl_{1} (L , h_{can})^{\wedge d} |_{X}
.
\]
It has the total volume $1$
and is \emph{smooth} and \emph{positive} on $X$.

We recall what 
the equidistribution theorem says.
Here we suppose that $K$ is a number field.
Let $(x_{l})_{l \in \NN}$ be a generic sequence of small points.
Let $\sigma$ be an archimedean place,
$X_{\sigma}$ the complex analytic space
of $X$ over $\sigma$,
and let $L_{\sigma}$ be the restriction of $L$ to $X_{\sigma}$.
Roughly speaking,
the equidistribution theorem says
that
the Galois orbits of $(x_{l})_{l}$,
approximatively as $l \to \infty$,
are
equidistributed in $X_{\sigma}$
with respect to
$\mu_{X_{\sigma},L_{\sigma}}$.

We can now recall the proof
in the arithmetic case due to Ullmo and Zhang.
The proof is done by contradiction.
Suppose we have a counterexample $X$ for the Bogomolov conjecture.
Then,
taking the quotient by the stabilizer of $X$ if necessary,
we can easily reduce ourselves to the case 
where
the stabilizer is trivial.
Further, we may assume $d := \dim X > 0$,
since the case of $d = 0$ is rather obvious.
For a large
$N \in \NN$, 
we can see that the morphism
\[
\alpha : X^{N} \to A^{N-1},
\quad
\alpha (x_{1}, \ldots , x_{N}) = 
(x_{2} - x_{1}, \ldots ,x_{N} - x_{N-1}).
\]
gives a birational morphism
$X^{N} \to \alpha (X^{N})$.
We fix such an $N$, writing $X ' := X^{N}$ and $Y := \alpha (X')$
for simplicity.
We
take dense Zariski-open subsets
$U \subset X'$ and $V \subset Y$
such that $\alpha$
induces an isomorphism between them.
Let $L'$ and $M$ be even ample line bundles on $X'$ and $Y$
respectively.
Then
we can see that $X'$ is again a counterexample
for the Bogomolov conjecture
with respect to the line bundle 
$L'$.
Therefore we can find a
generic sequence of small points
$(x_{l})_{l \in \NN}$,
and we may assume they sit in $U$.
Moreover, 
we can
see that
the image $(\alpha (x_{l}))_{l \in \NN}$ is also
a generic sequence of small points.
By virtue of the equidistribution theorem,
$(x_{l})_{l \in \NN}$ and $(\alpha (x_{l}))_{l \in \NN}$
are equidistributed in $X'$ and $Y$
with respect to
$\mu_{X'_{\sigma} ,L'_{\sigma}}$
and
$\mu_{Y_{\sigma} , M_{\sigma}}$
respectively,
for an archimedean place $\sigma$.
Furthermore since
$\alpha$ gives an isomorphism between $U$ and $V$,
we can conclude
\[
\mu_{X'_{\sigma} ,L'_{\sigma}} |_{U}
=
\alpha^{\ast}
(\mu_{Y_{\sigma} , M_{\sigma}}|_{V})
.
\]
Since both
$\mu_{X'_{\sigma} ,L'_{\sigma}}$
and
$\alpha^{\ast}
(\mu_{Y_{\sigma} , M_{\sigma}})$
are smooth forms,
we have
\[
\mu_{X'_{\sigma} ,L'_{\sigma}} 
=
\alpha^{\ast}
(\mu_{Y_{\sigma} , M_{\sigma}})
\]
on $X_{\sigma}$.
The right-hand side however cannot be positive over the
diagonal of $X' = X^{N}$.
This is
a contradiction since 
the left-hand side is positive.

How about the case of Gubler?
In contrast to the arithmetic case,
there are no archimedean places
in the geometric case.
That 
fact had prevented us from
enjoying an analogous proof
of the arithmetic case.
To overcome that difficulty,
Gubler used
non-archimedean analytic spaces over
a non-archimedean place
and their tropicalizations.

Let $X \subset A$ be a closed subvariety
of dimension $d$.
To a place
$v$ of $K$,
it is well known that
the Berkovich analytic spaces 
$X_{v} \subset A_{v}$ can be associated.
Gubler defined the canonical 
Chambert-Loir measure
$\mu_{X_{v}, L_{v}}$
on $X_{v}$.
Suppose here that $A_{v}$ is totally degenerate.
Then Gubler defined the
tropicalization
$X_{v}^{\trop}$
of $X_{v}$, which is denoted by $\overline{\val} (X_{v})$
in his article,
and
showed that
it is
a ``$d$-dimensional polytope''.
This plays the role of a counterpart of the complex space
over an archimedean place.
Furthermore he investigated in detail
the push-out $\mu_{X_{v}, L_{v}}^{\trop}$ to the tropicalization
of $\mu_{X_{v}, L_{v}}$,
describing it very concretely.
In fact he showed that
it is 
a $d$-dimensional positive Lebesgue measure 
on
the equi-$d$-dimensional polytope
$X_{v}^{\trop}$.

Now the idea of Ullmo and Zhang can be applied to this situation.
If there is a counterexample to the Bogomolov conjecture,
we can make a construction similar to that of the arithmetic case;
there is a morphism
$\alpha : X' \to Y$,
where $X'$ and $Y$ are 
some closed subvarieties of abelian varieties,
such that
$X'$ is again a counterexample of dimension $d' > 0$
and that 
$\alpha$ is
a generically finite morphism and 
the image of the diagonal
by
$\alpha$ is one point.
There is a generic net of small points since $X'$
is a counterexample to the Bogomolov conjecture.
Tropicalizing them, we have
\[
\alpha^{\trop} : (X'_{v})^{\trop}
\to Y_{v}^{\trop}
,
\]
which is a morphism of polytopes.
Using the tropical equidistribution theorem
of Gubler
to
a generic net of small points,
we can obtain
\addtocounter{Claim}{1}
\begin{align} \label{intequi}
\alpha^{\trop}_{\ast} 
(\mu_{X_{v}, L'_{v}}^{\trop}) = \mu_{Y_{v}, M_{v}}^{\trop}
\end{align}
as well,
where $L'$ and $M$ respectively are
even ample line bundles as before.
On the other hand,
there is a $d'$-dimensional face
$E$ such that $F := \alpha^{\trop} (E)$
is a lower dimensional face
since the subset corresponding to the diagonal 
contracts to a point.
It is impossible.
In fact,
the left-hand side of (\ref{intequi})
has a positive measure at a lower dimensional $F$,
but the right one is the $d'$-dimensional usual Lebesgue measure
as mentioned above.
Thus a contradiction comes out.

It is natural to ask whether
or not the same strategy works well
in the non-totally degenerate case.
It is known that 
the canonical measure
$\mu_{X_{v}, L}$ 
exists on the analytic space $X_{v}$.
Gubler defined in \cite{gubler3}
the tropicalization $X_{v}^{\trop}$
and studied the push-out
$\mu_{X_{v}, L}^{\trop}$ 
of the canonical measure.
He proved that
$X_{v}^{\trop}$ has the structure of a simplicial set
and that $\mu_{X_{v}, L}^{\trop}$ can be described
as
\addtocounter{Claim}{1}
\begin{align} \label{intro1}
\mu_{X_{v} , L}^{\trop}
=
\sum_{i = 1}^{N} r_{i} \delta_{\Delta_{i}}
,
\end{align}
where $\Delta_{i}$ runs through faces
and $\delta_{\Delta_{i}}$ is a usual relative Lebesgue measure on
the simplex $\Delta_{i}$.
On the other hand,
he also proved in
\cite{gubler4}
the
equidistribution theorem 
which holds true in this situation.

Thus we seem to have everything we need for the Bogomolov conjecture,
but we do not in fact.
When we obtained the contradiction
by using the equidistribution theorem,
it was crucial 
that
the canonical form or the canonical measure 
is a ``regular'' one.
Indeed,
if the canonical form were not smooth 
or positive in the arithmetic case,
a contradiction would not come out.
In Gubler's case as well,
it was the key that
the tropicalization of the canonical measure
is the Lebesgue measure on the equi-$d'$-dimensional polytope.
In the general case however,
lower dimensional $\Delta_{i}$'s
often appear
in (\ref{intro1}),
and that is troublesome.
It is true that we can make the same situation as before, that is,
we have a morphism
$
\alpha^{\trop} :
(X'_{v})^{\trop} \to 
Y_{v}^{\trop}
$
and
$\alpha^{\trop}_{\ast} 
\left(
\mu_{X_{v}', L'_{v}}^{\trop} 
\right)
= \mu_{Y_{v}, M_{v}}^{\trop}$
if we have a counterexample,
but
it is not sufficient to reach a contradiction
because $\mu_{Y_{v}, M_{v}}^{\trop}$ may contain a relative
Lebesgue measure
with a lower dimensional support.

The new idea in
this paper
to avoid this difficulty is to focus on 
the minimal dimension of a component of
the support of
$\mu_{Y_{v}, M_{v}}^{\trop}$.
In fact, we will show that it is bounded below by
the
abelian part of $Y_{v}$
(cf. \S~\ref{abelianpart}).
Then, we will see that
the equidistribution method
works quite well under the condition of 
Theorem~\ref{thoremint}.


\subsection{Further argument}

This paper contains
an
appendix
due to W. Gubler.
In communicating with the author
on the first version of this paper,
he found a proof of the fact that
the minimal dimension of the support of the components
of $\mu_{X_{v} , L}^{\trop}$ for ample $L$
is exactly $\dim X - b (X_{v})$.
Although we do not need this detailed information
in the proof of Theorem~\ref{thoremint},
it is quite interesting and will play an
important role
for applying the
canonical measures.
In fact, we will
apply Corollary~\ref{tropical canonical measure}
to make a contribution to the geometric Bogomolov conjecture
as follows:

\begin{Theorem}[cf. Corollary~\ref{GBC:bleq1}] \label{GBC:bleq1int}
Let $A$ be an abelian variety.
Suppose that there exists a place 
$v$ such that 
$b (A_{v}) \leq 1$.
Then the geometric Bogomolov conjecture holds for $A$.
\end{Theorem}

It is needless to say that this theorem
in the case of $b (A_{v}) = 0$
is Gubler's theorem.


\subsection{Organization}

This article is organized as follows.
We will give some remarks on the trace of an abelian variety
in \S~1.
Those who are familiar with the trace
will not have to read this section.
In \S~2,
we will formulate the
geometric Bogomolov conjecture for
abelian varieties.
In \S~3, 
we will deduce some results concerning our conjecture for a  
curve $X$ from the known jacobian cases.
We will describe in \S~4 some basic properties of
Berkovich analytic spaces and their tropicalizations.
We will also note some properties of the canonical measures.
In \S~5, we will give the proof of our main result.
The appendix
due to Gubler
is put at the
last part of this paper.
A result there will be used
in \S~\ref{furtherresults}.

\subsection*{Acknowledgments}
The author would like to thank Professor Fumiharu Kato
for giving
lectures on non-archimedean geometry
to the author,
Professor Walter Gubler
for giving a lot of valuable comments
and encouragement
and for his appendix,
The author would like to also thank
the referee for valuable comments.
This work was partially supported by
KAKENHI(21740012).

\subsection{Conventions and terminology} \label{CT}

Let $k$ be an algebraically closed field,
$\mathfrak{B}$ an irreducible normal projective variety over $k$,
and let $\mathcal{H}$ be an ample line bundle on $\mathfrak{B}$.\footnote{
We assume $\mathfrak{B}$ to be a curve in \S~\ref{knownresults}.}
Let
$K$ be the function field 
of 
$\mathfrak{B}$,
and let $\overline{K}$ be an algebraic closure of $K$.
All of them are fixed throughout this paper.

For a finite extension $K'$ of $K$,
let $\mathfrak{B}_{K'}$ denote the normalization of $\mathfrak{B}$ in $K'$.
Let $M_{K'}$ denote the set of points in $\mathfrak{B}_{K'}$ of codimension one.
For any $w \in M_{K'}$, the local ring $\OO_{\mathfrak{B}_{K'}, w}$ 
is a discrete valuation
ring with the fraction field $K'$,
and the order function $\ord_{w} : (K')^{\times} \to \ZZ$ gives an 
additive discrete valuation.
If $K''$ is a finite extension of $K'$, then
we have a canonical finite surjective morphism 
$\mathfrak{B}_{K''} \to \mathfrak{B}_{K'}$,
which induces a surjective map $M_{K''} \to M_{K'}$.
Thus we have an inverse system
$\left( M_{K'} \right)_{K'}$, where $K'$ runs through the finite extensions
of $K$ in $\overline{K}$,
and hence
we define $M_{\overline{K}} := \varprojlim_{K'} M_{K'}$.
We call an element of $M_{\overline{K}}$ a \emph{place} of $\overline{K}$.
Each place $v = ( v_{K'} )_{K'} \in M_{\overline{K}}$  determines a unique
non-archimedean 
multiplicative value $| \cdot |_{v}$ on $\overline{K}$
in such a way that
the following conditions are satisfied.
\begin{itemize}
\item
The restriction of $| \cdot |_{v}$
to $K'$ is equivalent to the valuation 
associated with the order function $\ord_{v_{K'}}$.
\item
For any $x \in K^{\times}$, $|x|_{v} = e^{ - \ord_{v_{K}} x}$.
\end{itemize}
Through this correspondence, we regard a place of $\overline{K}$ as
a valuation of $\overline{K}$.
For a $v \in M_{\overline{K}}$,
let $\overline{K}_{v}$ denote the completion of $\overline{K}$
with respect to $v$.
It is an algebraically closed field complete with respect to
the non-archimedean valuation $| \cdot |_{v}$.

For each $v_{K} \in M_{K}$,
let $| \cdot |_{v_{K} , \mathcal{H}}$ be the valuation normalized in such a way that
\[
| x |_{v_{K} , \mathcal{H}} := e^{- (\ord_{v_{K}} x)( \deg_{\mathcal{H}} v_{K}) }
,
\]
where $\deg_{\mathcal{H}} v_{K}$
stands for the degree of the closure of $v_{K}$
in $\mathfrak{B}$ with respect to $\mathcal{H}$.
It is well known that the set 
$\mathfrak{V} :=
\{ | \cdot |_{v_{K} , \mathcal{H}} \}_{v \in M_{K}}$ 
of valuations satisfies the product formula,
and hence we can define the notion of heights
with respect to this set of valuations,
namely, an absolute logarithmic height
with respect to $\mathfrak{V}$
(cf. \cite[Chapter~3 \S~3]{lang2}).
The ``height'' in this article
always means this height.

Let $F / k$ be any field extension.
For a scheme $X$ over $k$, we write 
$X_{F} := X \times_{\Spec k} \Spec F$.
If $\phi : X \to Y$ is a morphism of schemes over $k$,
we write $\phi_{F} : X_{F} \to Y_{F}$ for the base extension
to $F$.


\section{Descent of the base field of abelian varieties}
\label{appendix}

Let $F / k$ be any field extension.
We discuss in this section when an abelian variety
over $F$ can be defined over $k$,
and give a remark on the trace of an abelian variety.

Let us begin with a lemma.
\begin{Lemma} \label{etaleisogeny}
Let $A$ and $B$ be abelian varieties over $F$
and $k$ respectively.
If $\phi : B_{F} \to A$
is an \'etale isogeny,
then $A$ 
and $\phi$ are defined over $k$:
precisely,
there
exists a subgroup scheme $G$ of $B$ over $k$ such that
$\Ker \phi = G_{F}$
and hence
$A = (B / G)_{F}$.
\end{Lemma}

\begin{Pf}
Let $N$ be the degree of $\phi$
and let
$(B_{F}) [N]$
be the kernel of the $N$-times homomorphism.
Then $\Ker \phi \subset (B_{F}) [N]_{\mathrm{red}}$
since $\Ker \phi$ is reduced by our assumption.
Taking account that
the field extension $F / k$ is regular,
we
have 
\[
\Ker \phi \subset (B_{F}) [N]_{\mathrm{red}}
=
((B [N])_{F})_{\mathrm{red}}
=
\left( (B [N])_{\mathrm{red}}
\right)_{F}
,
\]
which tells us that
$\Ker \phi$
is defined over $k$,
namely,
there exists a subgroup scheme $G$ of $B$ over $k$ such that
$\Ker \phi = G_{F}$.
\end{Pf}

We recall here a 
quite fundamental theorem due to Chow:

\begin{citeTheorem}[cf. 
II \S~1 Theorem~5 in \cite{lang1}] \label{chow'stheorem}
Let $A$ be an abelian variety over $k$
and let $B$ be an abelian subvariety of $A_{F}$.
Then there exists an abelian subvariety $B' \subset A$
with $B'_{F} = B$.
\end{citeTheorem}

We can now show the following 
slight generalization of
Theorem~\ref{chow'stheorem}:

\begin{Proposition} 
\label{generalizedChow}
Let $A$ be an abelian variety over $k$
and let $G$ be a 
reduced closed subgroup of $A_{F}$.
Then there exists a
closed subgroup $G'$ of $A$
with $(G')_{F} = G$.
\end{Proposition}

\begin{Pf}
By Theorem~\ref{chow'stheorem},
there exists an abelian subvariety 
$G^{\circ} \subset A$ such that
$G^{\circ}_{F}$ is the identity component of $G$.
Consider the natural homomorphism
$\phi' : (A / G^{\circ})_{F} \to (A_{F}) / G$.
It is an \'etale isogeny
since $G$ is reduced, so by Lemma~\ref{etaleisogeny},
there exist an abelian variety $H$ over $k$
and a homomorphism $\psi : A / G^{\circ} \to H$
such that $\psi_{F}$ coincides with
$\phi'$.
Now let $G'$ be the kernel of 
the composition
$
A \to A / G^{\circ} \to H$.
Then we immediately find $G'_{F} = G$.
\end{Pf}

Let $B$ be an abelian variety over $k$
and let $\phi : B_{F} \to A$ be a smooth homomorphism
between abelian varieties over $F$.
Then,
as a corollary of Proposition~\ref{generalizedChow},
we can take an abelian variety $A'$ over $k$
and a homomorphism $\phi' : B \to A'$ such
that $A = A'_{F}$ and $\phi '_{F} = \phi$.
In fact,
there exists a
reduced closed subgroup
$G'$ of $B$ with $\Ker \phi = G'_{F}$
by Proposition~\ref{generalizedChow}.
Then 
$A ' := B / G'$
suffices our requirement.

Next we give a remark on the Chow trace.
Let $A$ be an abelian variety over $F$.
Recall that a
pair 
$\left( A^{F / k}, \Tr_{A}^{F/k} \right)$ of an abelian variety
$A^{F / k}$ over $k$ and a homomorphism 
$\Tr_{A}^{F/k} : (A^{F / k})_{F} \to A$
over $F$
is called a \emph{$F/k$-trace},
or \emph{Chow trace}, if it satisfies the following
universal property:
for any abelian variety $B$ over $k$
and for any homomorphism $\phi : B_{F} \to A$,
there exists a unique homomorphism
$\phi' : B \to A^{F / k}$ over $k$ such that
$\Tr_{A}^{F/k} \circ 
\phi'_{F} = \phi$
(cf. \cite {lang1} and \cite{lang2}).

\begin{Lemma} \label{finiteinsep}
$\Tr_{A}^{F/k}$ is finite and purely inseparable.
\end{Lemma}

\begin{Pf}
By virtue of Proposition~\ref{generalizedChow},
we can take
a closed subgroup $G' \subset A^{F/k}$ such that
$G'_{F} = 
\left( \Ker \Tr_{A}^{F/k} \right)_{\mathrm{red}}$.
Let
$\pi : A^{F/k} \to A^{F/k} / G' =: B$ be the quotient by $G'$.
Then we have naturally a homomorphism
$\phi : B_{F} 
\to A$.
By the universal property,
we obtain the factorization $\phi' : B \to A^{F/k}$
over $k$,
and the universality also
says that $\phi' \circ \pi  = \mathrm{id}_{A^{F/k}}$.
That concludes $\pi$ is an isomorphism and hence
$\left( \Ker \Tr_{A}^{F/k} \right)_{\mathrm{red}} = 0$,
namely,
$\Tr_{A}^{F/k}$ is finite and purely inseparable.
\end{Pf}

The uniqueness of the $F/k$-trace
is immediate from the definition.
We can find 
in \cite{lang1}
a proof for
the existence,
but we should note one thing:
in the definition of $F/k$-trace
of
\cite[VIII \S 8]{lang1},
Lang assumed that $\Tr_{A}^{F/k}$ is finite.
This assumption is not necessary
since it follows from the definition automatically
by virtue of Lemma~\ref{finiteinsep}
(cf. \cite[the last line in p.138]{lang2}).

Let $\phi : A \to B$ be a homomorphism of abelian varieties over $F$.
By the universality of the trace,
we see that
$\phi$ induces a homomorphism
$\Tr^{F/k} ( \phi ) : A^{F/k} \to B^{F/k}$.

\begin{Lemma} \label{surjectivity-trace}
With the notation above,
we
assume that $\phi$ is surjective.
Then $\Tr^{F/k} ( \phi ) : A^{F/k} \to B^{F/k}$ is surjective.
\end{Lemma}

\begin{Pf}
We can take an abelian subvariety $A_{0} \subset A$ such that
$\phi |_{A_{0}} : A_{0} \to B$ is an isogeny 
(cf. \cite[II, \S~1, Theorem~6]{lang1}).
Since $\Tr^{F/k} ( \phi |_{A_{0}} ) : A_{0}^{F/k} \to B^{F/k}$ 
factors through $\Tr^{F/k} ( \phi )$,
it is enough to show that $\Tr^{F/k} ( \phi |_{A_{0}} )$ is surjective.
Therefore,
we may assume that $\phi$ is an isogeny,
and hence $A = A_{0}$, at the beginning.
Then we can
take an isogeny $\psi : B \to A$,
and consider an isogeny $\theta := \phi \circ \psi : B \to B$.
Since $\theta \circ \Tr_{B}^{F/k} = \Tr_{B}^{F/k} \circ \Tr^{F/k} ( \theta)$
and since $\theta \circ \Tr_{B}^{F/k}$ is finite,
we see that
$\Tr^{F/k} ( \theta) : B^{F/k} \to B^{F/k}$ is finite.
Therefore
it is surjective,
and since
$\Tr^{F/k} ( \theta) = \Tr^{F/k} ( \phi ) \circ \Tr^{F/k} (\psi)$,
we conclude that $\Tr^{F/k} ( \phi )$ is surjective.
\end{Pf}

\section{Geometric Bogomolov conjecture}
In this section, we give a precise
formulation of the geometric
Bogomolov conjecture for abelian varieties.
\subsection{Small points} \label{set-up}
Let $A$ be an abelian variety over $\overline{K}$.
For an even ample line bundle $L$ on $A$,
let us consider the canonical height function $\hat{h}_{L}$.
It
is known to be a 
semi-positive quadratic form on $A (\overline{K})$.
Let $X$ be a closed subvariety of $A$.
For each $\epsilon > 0$,
we put
\[
X (\epsilon ; L)
:=
\left\{
x \in X (\overline{K})
\left|
\hat{h}_{L} (x) \leq \epsilon
\right.
\right\}
.
\]

\begin{Lemma} \label{image-counterexample}
Let $A$ and $X$ be as above.
Let $\phi : A \to B$ be a 
homomorphism
of abelian varieties over $\overline{K}$
and 
put $Y := \phi (X)$.
Let $L$ and $M$ be even ample line bundles on
$A$ and $B$ respectively.
Then if
$X (\epsilon ; L)$ is Zariski-dense in $X$
for any $\epsilon > 0$,
then
$Y (\epsilon' ; M)$ is Zariski dense in
$Y$
for any $\epsilon' >0$.
\end{Lemma}

\begin{Pf}
Since $L$ is ample,
we can take a positive integer $n$ such that
$L^{\otimes n} \otimes \phi^{\ast} (M)^{-1}$
is ample.
Then we have
$n \hat{h}_{L} \geq \phi^{\ast} \hat{h}_{M}$,
and hence
\begin{multline*}
Y (\epsilon ; M) 
=
\phi
\left(
\left\{ x \in X(\overline{K})
\left|
\hat{h}_{M} ( \phi (x)) \leq \epsilon
\right.
\right\}
\right)
\\
\supset
\phi
\left(
\left\{ x \in X(\overline{K})
\left|
n \hat{h}_{L} (x) \leq \epsilon
\right.
\right\}
\right)
=
\phi (X (\epsilon / n ; L))
.
\end{multline*}
The right-hand side is Zariski dense
in $Y$
by our assumption,
which leads us to our assertion.
\end{Pf}

Let $L_{1}$ and $L_{2}$ be even ample line bundles on $A$.
Then 
$X (\epsilon ; L_{1})$ is Zariski dense
for any $\epsilon > 0$ if and only if
so is $X (\epsilon' ; L_{2})$ for any $\epsilon' > 0$,
by virtue of the above lemma.
Accordingly, 
the following definition makes sense:
\begin{Definition}
We say that
\emph{X has dense small points}
if $X (\epsilon ; L)$ is Zariski dense
in $X (\overline{K})$ for any $\epsilon > 0$ and 
for some, hence any, even ample line bundle $L$ on $A$.
\end{Definition}

We end this subsection with
the following two basic lemmas
on small points,
which will be used later:

\begin{Lemma} \label{finite-homomorphism}
Let $\phi : A \to B$ be a 
isogeny
of abelian varieties over $\overline{K}$.
Let $X \subset A$
be a closed subvariety
and
put $Y := \phi (X)$.
Then 
$X$ has dense small points if and only if $Y$
has dense small points.
\end{Lemma}

\begin{Pf}
The ``only if'' part is immediate from Lemma~\ref{image-counterexample}.
Let us show the ``if'' part.
Let $M$ be an even ample line bundle on $B$.
Then $L := \phi^{\ast} M$ is also even and ample and we have
$
\phi (X (\epsilon ; L)) = Y (\epsilon ; M)$.
Then if $X (\epsilon ; L)$ is not Zariski dense for any $\epsilon > 0$,
then neither is not $Y (\epsilon ; M)$
since $\phi$ is finite.
This proves the ``if'' part,
and our lemma.
\end{Pf}

\begin{Lemma} \label{product-dense}
Let $A$ and $B$ be abelian varieties over $\overline{K}$
and let $X \subset A$ and 
$Y \subset B$
be closed subvarieties.
If $X$ and $Y$ have dense small points,
then the closed subvariety 
$X \times Y \subset A \times B$
also has dense small points.
\end{Lemma}

\begin{Pf}
Let $p : A \times B \to B$ and $q : A \times B \to B$
be the canonical projections.
For even ample line bundles $L$ and $M$ on $A$
and $B$ respectively, 
we write $L \boxtimes M := p^{\ast}{L} \otimes q^{\ast}{M}$.
It
is even ample and 
we have $\hat{h}_{L \boxtimes M} 
= p^{\ast} \hat{h}_{L} + q^{\ast} \hat{h}_{M}$.
Accordingly we have
\[
(X \times Y) (2 \epsilon ; L \boxtimes M)
\supset
X (\epsilon ;{L}) \times Y (\epsilon ;{M})
,
\]
and hence we obtain our assertion.
\end{Pf}

\subsection{Special subvarieties and the conjecture}
\label{specialsub}

For an abelian variety $A$ 
over $\overline{K}$,
let $\left( A^{\overline{K}/k} ,
\Tr_{A}^{\overline{K}/k} \right)$ denote the 
$\overline{K}/k$-trace of $A$
as in \S~\ref{appendix}.
Since $A^{\overline{K}/k}$ is defined over $k$,
we can consider $k$-points.
We note $A^{\overline{K}/k} (k) \subset 
A^{\overline{K}/k} \left( \overline{K} \right)$.

\begin{Definition} \label{definition-specialsubvariety}
Let $A$ be an abelian variety over $\overline{K}$
and let
$X \subset A$ be an
irreducible closed subvariety. 
Let $G_{X}$ be the stabilizer of $X$.
We call $X$ a \emph{special subvariety}
if there exist a torsion point $\tau \in A (\overline{K})_{tors}$,
and a closed subvariety
$X' \subset A^{\overline{K}/k}$
over $k$
such that
\[
X =
G_{X} +
\Tr_{A}^{\overline{K}/k}
\left(
{X'}_{\overline{K}}
\right)
+ \tau
,
\]
where $G_{X}$ is the stabilizer of $X$.
A point $x \in A ( \overline{K} )$ is called a 
\emph{special point} of $A$ if $\{ x \}$ is a special subvariety of
$A$.
\end{Definition}

Let $A_{sp}$ denote the set of special point of $A$.
It immediately follows from the definition that
\addtocounter{Claim}{1}
\begin{align} \label{specialnoshiki}
A_{sp} = A (\overline{K})_{tors} +
\Tr_{A}^{\overline{K}/k} \left(
A^{\overline{K}/k} (k)
\right)
.
\end{align}
Further,
let $L$ be an even ample line bundle.
Then we have
\addtocounter{Claim}{1}
\begin{align} \label{heightzeropoint}
A_{sp} 
= 
\left\{ x \in A (\overline{K}) 
\left|
\hat{h}_{L} (x) = 0
\right.
\right\}
.
\end{align}
by
\cite[Theorem 4.5 and 5.4.2]{lang2}.

\begin{Remark} \label{saigonitsuketaremark}
Let $X$ be an irreducible closed subvariety of $A$ and
let $\sigma \in A ( \overline{K})$ be a special point.
Then
$X$ is a special subvariety if and only if $X + \sigma$ is a
special subvariety,
which follows from the definition and (\ref{specialnoshiki}).
\end{Remark}


The following assertion says that
the
special subvarieties
have a dense subset of points of height zero:

\begin{Proposition} \label{dinseness-specialpoints}
If $X$ is a special 
subvariety of $A$,
then
$X \left( \overline{K} \right) \cap A_{sp}$ is dense in $X$.
\end{Proposition}

\begin{Pf}
There exists a closed subvariety $X' \subset A^{\overline{K}/k}$ and a torsion
point $\tau \in A( \overline{K})$ such that
\[
X = G_{X} + \Tr_{A}^{\overline{K}/k} \left( X'_{\overline{K}} \right) + \tau
.
\]
Then 
\[
G_{X} ( \overline{K})_{tors} + \Tr_{A}^{\overline{K}/k} \left( X'_{\overline{K}} (k)
\right) + \tau
\]
is a dense subset of $X$ and contained in $A_{sp}$.
\end{Pf}

In particular, we have the following:

\begin{Corollary} \label{special-small}
A special subvariety has dense small points.
\end{Corollary}

Now let us propose the statement of our 
geometric Bogomolov conjecture for abelian varieties,
which insists that the converse of Corollary~\ref{special-small}
should hold true:

\begin{Conjecture} [Geometric Bogomolov conjecture
for abelian varieties]
\label{GBC}
Let $A$ be an abelian variety over $\overline{K}$
and let $X$ be an irreducible closed subvariety of $A$.
Then,
$X$ should not have dense small points unless
it is a special subvariety.
\end{Conjecture}

We will see in Lemma~\ref{zerodimensional}
that
the above conjecture is easily verified
when $\dim X / G_{X} = 0$.
In the case where $\dim X / G_{X} > 0$,
the conjecture is not trivial at all.

Let us show some properties
on special points and special varieties
which will be needed.

\begin{Lemma} \label{surjectivity-specialpoints}
Let $\phi : A \to B$ be a surjective homomorphism
of abelian varieties over $\overline{K}$.
Then it induces a surjective homomorphism
$A_{sp} \to B_{sp}$.
\end{Lemma}

\begin{Pf}
We can show
$\phi \left( A \left( \overline{K} \right)_{tors} \right)
= B \left( \overline{K} \right)_{tors}$.
In fact, ``$\subset$" is trivial.
To show the other inclusion,
let us take an abelian subvariety $A_{0} \subset A$ such
that $\phi|_{A_{0}} : A_{0} \to B$ is finite and surjective.
Then we have
$\phi ( A_{0} ( \overline{K})_{tors} ) = B (\overline{K})_{tors}$,
which implies
$\phi \left( A \left( \overline{K} \right)_{tors} \right)
= B \left( \overline{K} \right)_{tors}$.
On the other hand, $\phi$ restricts to a surjective morphism
\[
 \Tr_{A}^{\overline{K}/k}
\left(A^{\overline{K}/k} (k) \right) \to 
\Tr_{B}^{\overline{K}/k} \left( B^{\overline{K}/k}  (k)
\right)
\]
by Lemma~\ref{surjectivity-trace}.
Accordingly, we obtain our lemma
by (\ref{specialnoshiki}).
\end{Pf}

\begin{Proposition} \label{specialsubvariety}
Let $X \subset A$ be a closed subvariety
and let $G_{X}$ be the stabilizer of $X$.
Put $B := A / G_{X}$
and $Y := X / G_{X}$.
Then the
following statements are equivalent to each other:
\begin{enumerate}
\renewcommand{\labelenumi}{(\alph{enumi})}
\item
$X$ is a special subvariety of $A$.
\item
$Y$ is a special subvariety of $B$.
\item
There exist an abelian variety $C$ over $k$,
a homomorphism 
$\phi : C_{\overline{K}} \to B$,
a closed subvariety $Z' \subset C$,
and
a special point $\sigma ' \in B ( \overline{K} )$ such that
$Y = \phi ({Z'}_{\overline{K}}) + \sigma '$.
\item
There 
exist a variety $W'$ over $k$,
a $k$-point $w_{0} \in W' (k)$,
a special point $\sigma \in Y(\overline{K})$
and a surjective morphism 
$\psi : {W'}_{\overline{K}} \to 
Y$
such that
$\psi (w_{0}) = \sigma$.
\end{enumerate}
\end{Proposition}

\begin{Pf}
We first show that (a) is equivalent to (b).
Let $\phi : A \to B$ be the canonical surjective homomorphism.
In order to show that (a) implies (b),
assume that $X$ is a special subvariety of $A$.
Then, we can 
take a closed subvariety $X' \subset A^{\overline{K}/k}$ and a torsion
point $\tau \in A ( \overline{K})_{tors}$ such that
$X = G_{X} + \Tr_{A}^{\overline{K}/k} \left( X'_{\overline{K}} \right) + \tau$.
Let $\Tr ( \phi ) : A^{\overline{K}/k} \to B^{\overline{K}/k}$
be the induced homomorphism
and
put $Y' :=
\Tr ( \phi ) ( X' )$.
Then, we have
\[
\phi
\left(
\Tr_{A}^{\overline{K}/k} \left( X'_{\overline{K}} \right)
\right)
=
\Tr_{B}^{\overline{K}/k} ( Y'_{\overline{K}} )
\]
and hence
\[
Y = \phi (X) = \phi
\left(
\Tr_{A}^{\overline{K}/k} \left( X'_{\overline{K}} \right)
\right)
+
\phi ( \tau )
=
\Tr_{B}^{\overline{K}/k} ( Y'_{\overline{K}} )
+
\phi ( \tau )
.
\]
Thus we find 
that $Y$ is a special subvariety of $B$.

In order to show that
(b) implies (a) next,
suppose that $Y$ is special.
Then, there exists a closed subvariety 
$Y'$ of $B^{\overline{K}/k}$ and a torsion point $\tau \in B ( \overline{K})$
such that $Y = \Tr_{B}^{\overline{K}/k}
\left(
Y'_{\overline{K}}
\right)
+ \tau$.
Since $\phi$ restricts to a surjective homomorphism 
$A ( \overline{K})_{tors} \to B ( \overline{K} )_{tors}$,
there exists a torsion point $\tau' \in A ( \overline{K})_{tors} $
with $\phi ( \tau' ) = \tau$.
Since $\Tr ( \phi )$ is surjective
by Lemma~\ref{surjectivity-trace},
we can take a closed subvariety $X' \subset A^{\overline{K}/k}$
with $\Tr ( \phi ) ( X' ) = Y'$.
Then, 
\[
\phi 
\left(
\Tr_{A}^{\overline{K}/k} \left( X'_{\overline{K}} \right)
+ \tau'
\right)
=
\Tr_{B}^{\overline{K}/k} ( Y'_{\overline{K}} ) + \tau
= Y
,
\]
and
therefore we find
\[
X = \phi^{-1} ( Y ) =
G_{X} +
\Tr_{A}^{\overline{K}/k} \left( X'_{\overline{K}} \right)
+ \tau'
.
\]
Thus we show that $X$ is a special subvariety of $A$.

The implication from (b) to (c) 
is trivial (cf. Remark~\ref{saigonitsuketaremark}).
It is also easy to see that (c) implies (d).
In fact,
with the notation in (c),
we put $W' := Z'$ and fix any point $w_{0} \in W' (k)$.
Further we
put $\sigma := \phi (w_{0}) + \sigma'$ and let 
$\psi : W'_{\overline{K}} \to Y$ be a morphism
defined by $\psi (w') = \phi (w') + \sigma'$.
Then they satisfies all the conditions in (d).
It only remains to
show that (d) implies (b).

Let $W'$, $w_{0}$, $\sigma$ and $\psi$
be those as in (d).
For a fixed $y \in B (\overline{K})$,
we define $T_{y} : B \to B$ by $T_{y} (x) = x + y$.
First note
that
we can write
$\sigma = \Tr_{B}^{\overline{K}/k} (t) + \tau$ with some
$t \in B^{\overline{K}/k} (k)$
and $\tau \in B(\overline{K})_{tors}$
by
(\ref{specialnoshiki}).
Then,
by considering $T_{- \tau} (Y)$ and $T_{- \tau} \circ \psi$
instead of $Y$ and $\psi$ respectively,
we may assume that 
$\sigma = \Tr_{B}^{\overline{K}/k} (t)$.
Further, taking an alteration of $W'$ if necessary,
we may and do assume that $W'$ is nonsingular.

Let us consider the albanese morphism
\[
\alpha'_{w_{0}}
:
W' \to \Alb (W')
\]
with respect to the base point $w_{0}$.
Then
\[
\alpha_{w_{0}} :=
(\alpha'_{w_{0}})_{\overline{K}}
:
{W'}_{\overline{K}} \to \Alb (W')_{\overline{K}}
=
\Alb ({W'}_{\overline{K}})
\]
is the albanese morphism of ${W'}_{\overline{K}}$
with respect to $w_{0}$.
By applying the universal property of 
$\alpha_{w_{0}}$
to the morphism $T_{-\sigma} \circ \psi$,
we obtain a homomorphism
\[
\phi : \Alb ({W'}_{\overline{K}})
\to B
\]
with $\phi \circ \alpha_{w_{0}} = T_{-\sigma} \circ \psi$.
Then by the universal property of the $\overline{K}/k$-trace,
$\phi$ factors through the $\overline{K}/k$-trace,
that is,
there is a homomorphism
$\phi' : \Alb (W') \to B^{\overline{K}/k}$
such that
\[
\Tr_{B}^{\overline{K}/k}
\circ ({\phi'}_{\overline{K}})
=
\phi
.
\]
We now consider a closed subvariety
$Y' := \phi' ( \alpha'_{w_{0}} (W')) + t$ of $A^{\overline{K} / k}$.
Then we have 
\begin{align*}
\Tr_{B}^{\overline{K}/k}
({Y'}_{\overline{K}}) 
=
\phi (\alpha_{w_{0}} ({W'}_{\overline{K}}))
+ \Tr_{B}^{\overline{K}/k} (t) 
= 
(T_{-\sigma} \circ \psi) ({W'}_{\overline{K}}) + 
\sigma
= 
(Y - \sigma) + \sigma
= Y 
\end{align*}
as required.
\end{Pf}

Now we show that
Conjecture~\ref{GBC}
is true if $\dim X / G_{X} = 0$.

\begin{Lemma} \label{zerodimensional}
Let $A$ be an abelian variety over $\overline{K}$
and let $X \subset A$ be a irreducible closed subvariety
such that $\dim X / G_{X} = 0$.
If $X$ is not a special subvariety, then
it does not have dense small points.
\end{Lemma}

\begin{Pf}
We can write $X / G_{X} = \{ \sigma \}$.
If $X$ has dense small points,
then so does $X / G_{X}$ by Lemma~\ref{image-counterexample}.
Therefore $\sigma$ is of height zero
and hence it is a special point by (\ref{heightzeropoint}).
That implies that $X / G_{X}$ is special, and so is $X$
by Proposition~\ref{specialsubvariety}.
\end{Pf}


\section{Some results for curves} \label{knownresults}

In this section,
we recall some known results concerning the
geometric Bogomolov conjecture
for jacobian varieties
and give remarks on their consequences.\footnote{We will not
use the results of this section
in the sequel.}
Throughout this section,
we assume that $K$ is a function field of a curve,
namely,
$\mathfrak{B}$ in \S~\ref{CT}
is a nonsingular curve.

Let $C$ be a curve over $\overline{K}$,
and let $J_{C}$ be the jacobian variety of $C$.
For each divisor on $C$ of degree $1$,
let
$j_{D} : C \to J_{C}$ be the embedding defined by
$j_{D} (x) = D - x$.
We note $j_{D} (x) + \sigma = j_{D + \sigma} (x)$
for each $\sigma \in J_{C}$.
The following assertion is an immediate consequence of
the theorem of Zhang and that of Cinkir.
We recall here that
a curve $C$ over $\overline{K}$ is \emph{isotrivial}
if it is a base extension 
of a curve over $k$
to $\overline{K}$.

\begin{Proposition} \label{jacobian}
Fix $c_{0} \in C \left( \overline{K} \right)$.
For each $\sigma \in J_{C} \left( \overline{K} \right)$,
we put $X_{c_{0},\sigma}^{\pm} := [\pm 1] (j_{c_{0}} (C) + \sigma)$,
where $[\pm 1]$ is the $\pm 1$-multiplication on 
$J_{C}$.
\begin{enumerate}
\item
Suppose that $C$ is isotrivial,
and let
$Z'$ be a curve over $k$ with an isomorphism
$\psi : {Z'}_{\overline{K}} \cong C$.
We assume further that $c_{0} \in \psi (Z'(k))$.
Then $X_{c_{0},\sigma}^{\pm}$ has dense small points 
if and only if $\sigma$
is a special point.
\item
Assume $\ch k = 0$.
If $C$ is non-isotrivial,
then
$X_{c_{0},\sigma}^{\pm}$ does not have dense small points.
\end{enumerate}
\end{Proposition}

\begin{Pf}
It is enough to consider
$X_{c_{0},\sigma} := X_{c_{0},\sigma}^{+}$ only.
Taking a finite extension of $K$ if necessary,
we may assume $C$ is 
a curve defined over $K$
with stable reduction at any place,
and $c_{0} \in C (K)$.
Then the assertion (2) is immediate from \cite[Theorem 2.12]{cinkir}
and \cite[Theorem 5.6]{zhang1}.

To see the assertion (1),
we first note that the admissible pairing $(\omega_{a}, \omega_{a})$
vanishes in this case (cf. \cite{zhang1}).
By virtue of \cite[Theorem 5.6]{zhang1},
we find that $X_{c_{0},\sigma}$ has dense small points
if and only if
the canonical height of the point 
corresponding to the divisor class
$(2 g - 2 )(c_{0} + \sigma) - \omega_{C}$
in the jacobian
vanishes.
Since 
$(2 g - 2 ) c_{0} - \omega_{C}$ is a special point of $J_{C}$
by our assumption,
that
is equivalent to $\sigma$ being special
in this case by (\ref{heightzeropoint}).
Thus we obtain our assertion.
\end{Pf}

Let us here prove a
technical lemma:

\begin{Lemma} \label{translation-lemma}
Let
$X$ be a closed subvariety of $A$,
and let 
$H \subset A$ 
be 
an abelian subvariety.
Suppose that there exists an $x_{0} \in A \left( \overline{K} \right)$
with $X - x_{0} \subset H$
and that
$X$ has dense small points.
Then there exists a special point
$\sigma$
of $A$
such that
$X - \sigma \subset H$.
Moreover,
$X - \sigma$ has dense small points.
\end{Lemma}

\begin{Pf}
The last statement follows form
\cite[Theorem 4.5 and 5.4.2]{lang2}
since $X$ has dense small points.
To complete the proof,
we may assume $H \subsetneq A$.
Let $\phi : A \to A/H$ be the quotient.
Since $X - x_{0} \subset H$,
we have $\phi (X) = \phi (x_{0})$.
Since $X$ has dense small points,
$\phi (x_{0})$ is a special point 
by Lemma~\ref{image-counterexample}.
By virtue of Lemma~\ref{surjectivity-specialpoints},
there exists $\sigma \in A_{sp}$ with 
$\phi (\sigma) = \phi (x_{0})$.
Then we have $\phi (X) = \phi (\sigma)$, that is,
$X - \sigma \subset H$.
\end{Pf}

Now we can show the following assertion,
which is a partial answer to
the geometric Bogomolov conjecture
when the closed subvariety $X$ is a curve.

\begin{Proposition}[$\ch k = 0$]
 \label{GBC-for-curve-with-simple-jacobian}
Let $X$ be an irreducible closed subvariety of $A$ of dimension $1$,
and let $\nu : Y \to X$ be the normalization.
Let $J_{Y}$ be the jacobian variety of $Y$.
Suppose 
that
$J_{Y}$ is simple.
Then $X$ does not have dense small points
unless it is a special subvariety.
\end{Proposition}

\begin{Pf}
For a fixed $y_{0} \in Y (\overline{K})$,
we
put $x_{0} := \nu (y_{0})$ and $X_{0} := X - x_{0}$.
Then $0 \in X_{0} (\overline{K})$ and we have naturally
$\nu_{0} : Y \to X_{0}$ with $\nu (y_{0}) = 0$,
by composing the translation by $-x_{0}$
to $\nu$.
Then,
by the universality of the albanese variety $J_{Y}$,
we have a homomorphism $\phi : J_{Y} \to A$ such that
the
diagram
\begin{align*}
\begin{CD}
Y @>{j_{y_{0}}}>> J_{Y} \\
@V{\nu_{0}}VV @VV{\phi}V \\
X_{0} @>{\subset}>> A 
\end{CD}
\end{align*}
is commutative.
We require an additional condition on $y_{0}$ 
in case that $Y$ is an isotrivial curve:
we can take a variety
$Y'$ 
over $k$ and an isomorphism
$\psi : {Y'}_{\overline{K}} \cong Y$,
and
our requirement is $y_{0} \in \psi (Y' (k))$.

Under the setting above,
we will show Proposition~\ref{GBC-for-curve-with-simple-jacobian} 
by contradiction.
Suppose that $X$ is not a special subvariety
but it has dense small points.
Let $H$ be the image of the homomorphism $\phi$.
Then by Lemma~\ref{translation-lemma},
there is a $\sigma \in A_{sp}$
such that $X_{1} := X - \sigma \subset H$,
and
moreover $X_{1}$ has dense small points.
We put $z := \sigma - x_{0}$.
Then
we have
$X_{1} = X_{0} - z$ 
and $z \in H$.
We take $w \in J_{Y}$ with $\phi (w) = z$
and
consider $Y_{1} := j_{y_{0}} ( Y ) - w$.
Note
$\phi (Y_{1}) = X_{1}$.
The homomorphism $\phi$ is finite
since $J_{Y}$ is simple by our assumption.
Therefore, 
we see that $Y_{1}$ has dense small points 
by Lemma~\ref{finite-homomorphism}.

Here we divide ourselves into two cases.
One is the case where
$Y$ is non-isotrivial.
Then $Y_{1}$ cannot have dense small points
by Proposition~\ref{jacobian} (2),
hence the contradiction immediately comes out.

Let us consider the other case,
namely, the case where $Y$ is isotrivial.
Since $Y_{1}$ has dense small points 
and 
$Y_{1} = j_{y_{0}} (Y) - w$,
we see that $w$ is a special point
by Proposition~\ref{jacobian} (1).
That says that $z = \phi (w)$ is a special point,
which implies $X_{1}$ a special subvariety
by Proposition~\ref{specialsubvariety}.
Accordingly
$X = X_{1} + \sigma$ is also a special subvariety
(cf. Remark~\ref{saigonitsuketaremark}),
which
contradicts our assumption.
Thus we have proved our assertion.
\end{Pf}

\section{Preliminaries} \label{preliminary}

We fix our conventions and terminology.
When we write $\KK$,
it is a field 
which is complete with respect to a
non-trivial
non-archimedean absolute value $| \cdot | : \KK^{\times} \to \RR$.
Our $\overline{K}_{v}$,
which is the completion of
$\overline{K}$ with respect to 
a valuation $v$ as in \S~\ref{CT}, is a typical example of 
$\KK$.
Let $\Gamma := \{ - \log |a| \mid a \in \KK^{\times} \}$ be the value group
of $\KK$.
We put 
$
\KK^{\circ} := 
\{
a \in \KK \mid | a | \leq 1
\}
$,
the ring of integers of $\KK$,
and put
$
\KK^{\circ \circ} := 
\{
a \in \KK \mid | a | < 1
\}
$,
the maximal ideal of the valuation ring $\KK^{\circ}$.
Further we put
$\tilde{\KK} := \KK^{\circ} / \KK^{\circ \circ}$.
For an admissible formal scheme\footnote{Any formal scheme is supposed to be
an admissible formal scheme in this article.} $\XXX$
(cf.
\cite{gubler1, gubler3}),
we write $\tilde{\XXX} := \XXX \times_{\Spf \KK^{\circ}} 
\Spec \tilde{\KK}$.
For a morphism $\varphi : \XXX \to \YYY$
of admissible formal schemes,
we write $\tilde{\varphi}$ for the induced morphism
between their special fibers.

\subsection{Berkovich analytic spaces} \label{bersp}
In the theory of non-archimedean analytic geometry,
there are several kinds of ``visualization''
or in other words, some kinds of spaces
that realize the theory of non-archimedean analytic geometry.
In this article,
we adopt the spaces introduced by Berkovich
which are called \emph{Berkovich analytic spaces}.
When we say an \emph{analytic space},
it always means a Berkovich analytic space in this article.
In this subsection, we recall some notions and properties on
analytic spaces associated to
admissible formal schemes of
algebraic varieties,
as far as we use later.
For details, we refer to
his original papers
\cite{berkovich1,berkovich2,berkovich3,berkovich4}.
We also refer
to
Gubler's expositions in his papers \cite{gubler1,gubler3},
which
would be good reviews to this theory.

Let $\XXX$ be an
admissible formal scheme over $\KK^{\circ}$.
Then we can associate an analytic space $\XXX^{\an}$,
called the \emph{generic fiber} of $\XXX$.
For a given analytic space $X$, an admissible formal scheme
with the generic fiber $X$ is called a
\emph{formal model} of $X$.
There is a reduction map
$
\mathrm{red}_{\XXX} :
\XXX^{\an} \to \tilde{\XXX}
$.
Let $Z$ an irreducible component of $\tilde{\XXX}$
with the generic point $\xi_{Z}$.
Then there is a unique point $\eta_{Z} \in \XXX^{\an}$
with $\mathrm{red}_{\XXX} (\eta_{Z}) = \xi_{Z}$.
Thus we can naturally regard the generic
point of each irreducible component of the special fiber
as a point of the generic fiber.

We can also associate an analytic space to an
algebraic variety.
Let $X^{\an}$
denote
the analytic space associate to an
algebraic variety $X$ over $\KK$.
We have
naturally
$X (\KK) \subset X^{\an}$.
We should give a remark on the relationship between
the analytic space associated to
an algebraic variety and that done to an admissible formal scheme.
Let $X$ be an algebraic scheme over $\KK$.
Let $\mathcal{X}$ be a model of $X$,
that is,
$\mathcal{X}$ is a scheme
flat and of finite type
over $\KK^{\circ}$ with the generic fiber $X$.
Let $\hat{\mathcal{X}}$ be the formal completion with respect to
a nontrivial principal open ideal of $\KK^{\circ}$.
Then 
it is an admissible formal scheme
and
$\hat{\mathcal{X}}^{\an}$ is an analytic subdomain of
$X^{\an}$.
Moreover
if $\mathcal{X}$ is proper over $\KK^{\circ}$,
then $\hat{\mathcal{X}}^{\an} = X^{\an}$.

Let $X$ be a proper algebraic variety over $\KK$
and let $Y$ be its closed subvariety.
Let 
$\XXX$ 
be a formal model of $X^{\an}$.
Then there is a unique admissible formal subscheme
$\YYY \subset \XXX$ with $\YYY^{\an} = Y^{\an}$.
We call this $\YYY$ the \emph{closure} of $Y$ in $\XXX$.

Finally we fix a notation.
Let $X$ be an algebraic scheme over $\overline{K}$,
and let $v$ be a place of $\overline{K}$.
Then we have a Berkovich analytic space
associated to $X \times_{\overline{K}} \Spec \overline{K}_{v}$.
We denote it by $X_{v}$.
It is a typical analytic space
that we will mainly deal with in the sequel.

\subsection{Raynaud extension} 
\label{RtM}
For simplicity,
we assume further that $\KK$ is algebraically closed
from here on.
We recall here some notions of the Raynaud extensions
as far as needed in the sequel.
See \cite[\S 1]{BL} and \cite[\S 4]{gubler3} for details.

Let $A$ be an abelian variety over $\KK$.
According to \cite[Theorem 1.1]{BL},
there exists 
a unique analytic subgroup
$A^{\circ} \subset A^{\an}$
with a formal model $\mathscr{A}^{\circ}$ 
having the following properties:
\begin{itemize}
\item
$\mathscr{A}^{\circ}$ is a formal group scheme and
$(\mathscr{A}^{\circ})^{\an} \cong A^{\circ}$ as group analytic spaces.
\item
There is 
a short exact sequence
\addtocounter{Claim}{1}
\begin{align} \label{qcptformalR}
\begin{CD}
1 @>>> \mathscr{T}^{\circ} @>>> \mathscr{A}^{\circ} @>>>
\mathscr{B} @>>> 0 ,
\end{CD}
\end{align}
where $\mathscr{T}^{\circ}$ is a formal torus and 
$\mathscr{B}$ is a formal abelian variety.
\end{itemize}
By virtue of \cite[Satz 1.1]{bosch},
we see that such an
$\AAA^{\circ}$ is unique,
and $\TTT^{\circ}$
and $\BBB$ are also uniquely determined.
Taking the generic fiber of (\ref{qcptformalR}),
we have an exact sequence
\begin{align*}
\begin{CD}
1 @>>> T^{\circ} @>>> A^{\circ} @>{(q^{\circ})^{\an}}>> B @>>> 0
\end{CD}
\end{align*}
of group spaces.
We call $\TTT^{\circ}$, $\AAA^{\circ}$ and $\BBB$ 
are the \emph{canonical formal models}
of $T^{\circ}$, $A^{\circ}$ and $B$ respectively.
 
Naturally $T^{\circ}$ is a quasicompact subgroup
of the analytic torus 
$T$,
and hence we can obtain the push-out of the above extension:
\addtocounter{Claim}{1}
\begin{align} \label{Raynaudext}
\begin{CD}
1 @>>> T @>>> E @>{q^{\an}}>> B @>>> 0 .
\end{CD}
\end{align}
The natural morphism $A^{\circ} \to E$
is an immersion of analytic groups,
of which image is denoted 
by $E^{\circ}$.
\cite[Theorem 1.2]{BL} says that
the homomorphism $T^{\circ} \hookrightarrow A^{\an}$ extends uniquely
to a homomorphism $T \to A^{\an}$ and hence to 
a homomorphism
$p^{\an} : E \to A^{\an}$.
This $p^{\an}$
is called the \emph{Raynaud extension}
of $A$.
Note that we have a commutative diagram
\addtocounter{Claim}{1}
\begin{align} \label{opensubgroups}
\begin{CD}
E^\circ @>{\subset}>> E \\
@VV{\cong}V @VV{p^{\an}}V \\
A^{\circ} @>{\subset}>> A^{\an}
\end{CD}
\end{align}
of group spaces.
It is known that
$p^{\an}$ is a surjective homomorphism and moreover
$M := \Ker p^{\an}$ is a lattice.
Thus $A^{\an}$ can be described as the quotient of
$E$ by a lattice $M$.

We recall the valuation maps next.
Taking account that
the transition functions of the
$T$-torsor (\ref{Raynaudext}) can be valued in $T^{\circ}$,
we can define a continuous map
\[
\val : E \to \RR^{n}
,
\]
as in \cite{BL},
where $n := \dim T$ is called the \emph{torus rank} of $A$.
In fact,
we can
take an
analytic subdomain $V \subset B$
and a trivialization
\addtocounter{Claim}{1}
\begin{align} \label{triv-1}
(q^{\an})^{-1} (V)
\cong V \times T 
\end{align}
such that its  restriction induces a trivialization
\begin{align*} 
((q^{\circ})^{\an})^{-1} (V)
\cong V \times T^{\circ} .
\end{align*}
Let us
consider the composite
$
r_{V} : (q^{\an})^{-1} (V) \cong V \times T 
\to T
$ 
of (\ref{triv-1})
with the second projection.
We see that
if $x \in (q^{\an})^{-1} (V)$,
then
$
\val (x)
=
(v (r_{V} (x)_{1})
,
\ldots
,
v (r_{V} (x)_{n})
)
$,
where $r_{V} (x)_{j}$
is the $j$-th coordinate of $r_{V} (x)$.
Moreover, the
lattice $M$ is 
mapped by $\val$ to a 
lattice $\Lambda \subset
\RR^{n}$ and 
we have a diagram
\addtocounter{Claim}{1}
\begin{align}
\label{tropicalization}
\begin{CD}
E @>{\val}>> \RR^{n} \\
@VVV @VVV \\
A^{\an} @>{\overline{\val}}>> \RR^{n} / \Lambda
\end{CD}
\end{align}
that is commutative.
From the construction 
of the map $\val$,
we can see
\[
E^{\circ} = T^{\circ} \cap E = \val^{-1} (\mathbf{0})
.
\]

\subsection{Mumford models} \label{RtM2}
Let $\CCC$ be a 
$\Lambda$-periodic $\Gamma$-rational
polytopal decomposition of $\RR^{n}$
(cf. \cite[\S 6.1]{gubler1}).
Taking the quotient by $\Lambda$,
we have a polytopal decomposition of $\RR^{n} / \Lambda$.
Gubler  constructed the
\emph{Mumford model}\footnote{In this article, 
we consider Mumford models associated 
with a 
$\Lambda$-periodic
polytopal decomposition which is
\emph{$\Gamma$-rational} only.}
\[
p = p_{\CCC} : \EEE \to \AAA
\]
associated to $\CCC$.
We also call $\AAA$ the Mumford model of $A$.
We refer to \cite[\S 4]{gubler3} for details,
and recall some properties that will be needed:
\begin{itemize}
\item
The surjection $q^{\an} : E \to B$ extends to
$q : \EEE \to \BBB$ uniquely.
If $\TTT$ denote the closure of $T$ in $\EEE$,
then 
\addtocounter{Claim}{1}
\begin{align} \label{morphismq}
q : \EEE \to \BBB 
\end{align}
is a fiber bundle
with the fiber $\TTT$.
\item
The lattice $M$ acts freely on $\EEE$
and
$\EEE / M = \AAA$.
In particular
$p$ is a local isomorphism. 
\item
If $\CCC'$ is a polytopal decomposition of 
$\RR^{n}$ finer than $\CCC$,
and if $\EEE' \to \AAA'$ is the Mumford model associated to $\CCC'$,
then there are natural morphisms $\EEE' \to \EEE$
and
$\AAA' \to \AAA$
extending the identity homomorphisms.
\end{itemize}

Let 
$\EEE \to \AAA$ be a Mumford model
of the Raynaud extension $p^{\an} : E \to A^{\an}$
and let
$\red_{\EEE} : E \to \tilde{\EEE}$ be the reduction map.
By virtue of \cite[Proposition 4.8]{gubler3},
we see that
$\red_{\EEE} (E^{\circ}) = \red_{\EEE} (\val^{-1} (\mathbf{0}))$
is an open subset 
of $\tilde{\EEE}$.
Let $\EEE^{\circ}$ be the open formal subscheme of $\tilde{\EEE}$
with
$\tilde{\EEE^{\circ}} = \red_{\EEE} ( E^{\circ} )$.
The group structure of $E^{\circ}$ endows
$\EEE^{\circ}$ with a group structure.
Since we find
\[
\Ker \left(
\tilde{q}|_{\tilde{\EEE}^{\circ}} : \tilde{\EEE}^{\circ} \to \tilde{\BBB}
\right)
=\tilde{\TTT}^{\circ}
\]
by \cite[Remark 4.9]{gubler3},
we obtain
an exact sequence
\addtocounter{Claim}{1}
\begin{align} \label{formalraynaud2}
\begin{CD}
1 @>>> \mathscr{T}^{\circ}
@>>> \EEE^{\circ} @>q^{\circ}>>
\BBB @>>> 0
\end{CD}
\end{align}
of formal group schemes.
By the uniqueness property of $\AAA^{\circ}$,
we then conclude 
that $\EEE^{\circ} \cong \AAA^{\circ}$
and (\ref{formalraynaud2}) coincides with
(\ref{qcptformalR}).
Using this isomorphism, we can define
$\AAA^{\circ} \to \AAA$ to be the composite
\[
\AAA^{\circ} \cong \EEE^{\circ} \subset \EEE \stackrel{p}{\to} \AAA
.
\]
This morphism
$\AAA^{\circ} \to \AAA$ 
is an open immersion since $p|_{\EEE^{\circ}}$ is an isomorphism,
and
it
is an extension of $A^{\circ} \subset A$ to their formal models.
Thus we obtain a 
commutative diagram
\addtocounter{Claim}{1}
\begin{align} \label{modelofopensubgroups}
\begin{CD}
\EEE^\circ @>{\subset}>> \EEE \\
@V{p^{\circ}}V{\cong}V @VV{p}V \\
\AAA^{\circ} @>{\subset}>> \AAA ,
\end{CD}
\end{align}
which is a formal model of a diagram (\ref{opensubgroups}).

\subsection{Tropicalization}
Let $X$ be
a closed subvariety of $A$.
Then the image $\overline{\val} (X^{\an})$ is 
known to be a closed subset
of $\RR^{n} / \Lambda$.
We put
\[
X^{\trop} : = \overline{\val} (X^{\an})
,
\]
calling it the \emph{tropicalization} of $X$.
It is well known that
$X^{\trop}$ has
the structure of a polytopal set (cf. \cite[Theorem~1.1]{gubler3}).

The following assertion will be used in the proof of 
our main result.

\begin{Lemma} \label{product-trop}
Let $A_{1}$ and $A_{2}$ be abelian varieties over $\KK$
and let $X_{1} \subset A_{1}$ and $X_{2} \subset A_{2}$
be closed subvarieties.
Then we have naturally
\[
(X_{1} \times X_{2})^{\trop}
=
{X_{1}}^{\trop} \times {X_{2}}^{\trop} 
.
\]
\end{Lemma}

\begin{Pf}
From the definition of the tropicalization,
we immediately see 
\[
{A_{1}}^{\trop} \times {A_{2}}^{\trop} 
=
\RR^{n_{1}} / \Lambda_{1} \times \RR^{n_{2}} / \Lambda_{2}
=
\RR^{n_{1} + n_{2}} / \Lambda_{1} \times \Lambda_{2}
= (A_{1} \times A_{2})^{\trop}
,
\]
where $n_{i}$ is the torus rank of $A_{i}$
and $\Lambda_{i}$ is the lattice 
as in (\ref{tropicalization}) for $A_{i}$.
Both
$(X_{1} \times X_{2})^{\trop}$
and
${X_{1}}^{\trop} \times {X_{2}}^{\trop}$ 
are subsets of the above real torus.
On the other hand,
we have a natural surjective map
$(X_{1} \times X_{2})^{\trop} 
\to {X_{1}}^{\trop} \times {X_{2}}^{\trop}$
associated to the natural surjective continuous map
\[
\left| (X_{1} \times X_{2})^{\an} \right|
\to
\left| (X_{1})^{\an}
\right| 
\times 
\left| (X_{2})^{\an} \right|
,
\]
where $| X^{\an} |$ stands for the underlying topological
space of a Berkovich analytic space $X^{\an}$.
Thus we conclude $(X_{1} \times X_{2})^{\trop} 
= {X_{1}}^{\trop} \times {X_{2}}^{\trop}$.
\end{Pf}


\subsection{The dimension of the abelian part  of a closed subvariety}
\label{abelianpart}
In this subsection,
let $A$ be an abelian variety over $\KK$
and let $X \subset A$
be
an irreducible closed subvariety.

\begin{Lemma} \label{welldef-ab}
For $i = 0,1$,
let $p_{i} : \mathscr{E}_{i} \to \AAA_{i}$ be a Mumford model of 
the Raynaud extension
of $A$
and let $q_{i} : \EEE_{i} \to \BBB$ be the
morphism as \textup{(\ref{morphismq})}.
Let
$\XXX_{i}$
be
the closure
of $X$ in $\AAA_{i}$
and let $\YYY_{i}$ be
a quasicompact open subscheme
of $p^{-1} (\XXX_{i})$
such that
$p (\YYY_{i}) = \XXX_{i}$.
Then we have
$\dim \tilde{q_{0}} (\tilde{\YYY_{0}}) =
\dim \tilde{q_{1}} (\tilde{\YYY_{1}})$. 
\end{Lemma}

\begin{Pf}
We can take a Mumford model $p : \EEE \to \AAA$ such that 
$\AAA$ dominates both
$\AAA_{0}$ and $\AAA_{1}$.
Let $q : \EEE  \to \BBB$ be
the morphism as (\ref{morphismq}).
We also have a dominant morphism
$\XXX \to \XXX_{i}$ for $i = 0,1$,
where $\XXX$ is the closure of $X$ in $\AAA$.
Set $\YYY'_{i} := \XXX \times_{\XXX_{i}} \YYY_{i}$.
Then $\YYY'_{i}$ is a quasicompact open formal subscheme
of $p^{-1} (\XXX)$ 
such that 
$\YYY'_{i} \to \XXX$ 
is surjective.
Moreover
$\tilde{\YYY'_{i}} \to \tilde{\YYY_{i}}$ is surjective,
and hence $\dim ( \tilde{q} ( \tilde{\YYY'_{i}})) = 
\dim ( \tilde{q_{i}} ( \tilde{\YYY_{i}}))$.
Accordingly,
by pulling everything
back to
$\AAA$, 
we may assume 
$\AAA_{1} = \AAA_{0} = \AAA$ and hence
$\XXX_{1} = \XXX_{0} = \XXX$,
$\EEE_{1} = \EEE_{0} =: \EEE$
and $\YYY_{0}, \YYY_{1} \subset p^{-1} (\XXX)$.

Let us fix an irreducible component $W$ of $\tilde{\XXX}$.
There are irreducible components
$Z_{0}$ and $Z_{1}$
of $\tilde{\YYY_{0}}$ and $\tilde{\YYY_{1}}$
lying over $W$.
Since $\AAA = \EEE / M$,
there exist $m \in M$ such that
$Z_{1} \cap ({Z_{0}} + m)$
is a Zariski dense open subset of both 
$Z_{1}$
and
$Z_{0} + m$.
Accordingly,
\[
\dim \tilde{q} (Z_{1})
=
\dim \tilde{q}
(Z_{1} \cap (Z_{0} + m))
=
\dim \tilde{q} 
(Z_{0} + m)
=
\dim \tilde{q} (Z_{0})
.
\]
Consequently,
the number
$
\dim \tilde{q} (Z)$,
where $Z$ is an irreducible component over $W$,
depends only on $W$.
If we write $\alpha (W)$ for this number,
we see, for each $i$, that
\[
\dim \tilde{q}(\tilde{\YYY_{i}})
=
\max_{W} \alpha (W)
,
\]
where $W$ runs through the irreducible components
of the quasicompact scheme $\tilde{\XXX}$.
Thus
our assertion follows.
\end{Pf}

By virtue of the above lemma, we
can make the following definition:

\begin{Definition}
For
an irreducible closed subvariety
$X$ of $A$,
we define 
$b (X)$ to be the
number $\dim \tilde{q_{0}} (\tilde{\YYY_{0}}) =
\dim \tilde{q_{1}} (\tilde{\YYY_{1}})$
in Lemma~\ref{welldef-ab}.
We call it 
\emph{the dimension of the abelian part of $X$}.
\end{Definition}

Note
$b (A) = \dim B$,
where $B$ is the abelian part of the Raynaud extension of $A$
(cf. (\ref{Raynaudext})),
and
$b (A) = 0$ if and only if $A$ is
totally degenerate.


\begin{Lemma} \label{quotientb}
Let $X \subset A$ be an irreducible closed subvariety
and let $G_{X}$ be the stabilizer of $X$.
Then $b (X) \geq b (X / G_{X})$.\footnote{We can 
actually show that an equality $b (X) = b ( X / G_{X}) + b (G_{X})$
holds,
though we do not use it in this article.}
\end{Lemma}

\begin{Pf} 
Let $\phi : A \to A / G_{X}$ be the quotient homomorphism.
Then $\phi$ lifts to a homomorphism between 
the Raynaud extensions of $A$ and $A/G_{X}$
by \cite[Theorem 1.2]{BL}.
Therefore,
if $\BBB$ and $\CCC$ are the formal abelian varieties
such that $\BBB^{\an}$ and $\CCC^{\an}$
are the abelian parts of the Raynaud extensions of
$A$ and $A/G_{X}$ respectively,
then we have an induced homomorphism $\BBB \to \CCC$.
Now
our assertion follows immediately from the definition
of $b (X)$.
\end{Pf}

\subsection{Chambert-Loir measures}

The purpose of this subsection is 
to give a remark on
the product of
Chambert-Loir measures.
We refer to \cite[\S 3]{gubler3}
for all the notions such as
admissible metric and Chambert-Loir measures.

Let
$X$ be a projective variety over $\overline{K}$.
Recall that $\overline{K}_{v}$
denote the completion of $\overline{K}$
with respect to  a place $v \in M_{\overline{K}}$
(cf. \S~\ref{CT}),
and
$X_{v}$ the analytic space
associated to a algebraic variety $X \times_{\Spec \overline{K}} 
\Spec \overline{K}_{v}$ 
(cf. \S~\ref{bersp}).
To an admissibly metrized line bundle
$\overline{L}$ on $X$
(cf. \cite[\S 3.5]{gubler3}),
we can associate 
a Borel measure 
\[
\mu_{X_{v},\overline{L}} 
:=
\frac{1}{
\deg_{L} X}
\cherncl_{1} (\overline{L})^{\wedge d}
\]
on $|X_{v}|$
(cf. \cite[Proposition~3.8]{gubler3}),
where we emphasize with $| \cdot |$ that
$|X_{v}|$ is the underlying topological space
of the Berkovich analytic space $X_{v}$.

The following formula is the one
mentioned in \cite[\S 2.8]{chambert-loir} essentially,
but we restate it
with a proof for readers' convenience.

\begin{Proposition} [\S 2.8 in \cite{chambert-loir}]
\label{productmeasure}
Let $X$ and $Y$ be projective varieties over $\overline{K}$
and let $\overline{L}$ and $\overline{M}$ be 
admissibly metrized line bundles
on $X$ and $Y$ respectively.
Let 
$p$ and $q$ be the canonical projections from $X \times Y$
to $X$ and $Y$ respectively,
and let
$r : 
| X_{v} \times Y_{v} | \to |X_{v}| \times |Y_{v}|$ 
be the canonical continuous
map
induced from the projections.
Then we have
\[
\mu_{X_{v}, \overline{L}} \times \mu_{Y_{v} , \overline{M}} = 
r_{\ast} 
\left(
\mu_{X_{v} \times Y_{v} , \overline{L} \boxtimes \overline{M}}
\right)
,
\]
where $\overline{L} \boxtimes \overline{M} = 
p^{\ast} \overline{L} \otimes q^{\ast} \overline{M}$
and
$\mu_{X_{v}, \overline{L}} \times \mu_{Y_{v} , \overline{M}}$ is the
product measure on $|X_{v}| \times |Y_{v}|$.
\end{Proposition}

\begin{Pf}
First let us consider the case where the
admissible metric on $L$ and $M$ are the formal metrics
arising from models $(\XXX,\LLL)$ and $(\YYY,\MMM)$ respectively
(cf. \cite[\S 3]{gubler1}).
Then $\overline{L} \boxtimes \overline{M}$
is the formally metrized line bundles arising from the model
$(\XXX \times \YYY,  \LLL \boxtimes \MMM)$.
By virtue of \cite[Proposition 3.11]{gubler1},
we have explicit formulas
\begin{align*}
\mu_{\overline{L}} &= 
\frac{1}{\deg_{L} X} \sum_{A \in \mathrm{Irr} 
\left(
\tilde{\XXX}
\right)} 
\left( \deg_{\LLL} A \right) \delta_{\eta_{A}} \\
\mu_{\overline{M}} &= \frac{1}{\deg_{M} Y} 
\sum_{B  \in \mathrm{Irr} 
\left(
\tilde{\YYY}
\right)} 
\left( \deg_{\MMM} B \right) \delta_{\eta_{B}} \\
\mu_{\overline{L} \boxtimes \overline{M}} 
&= \frac{1}{\deg_{L \boxtimes M} X \times Y} 
\sum_{C  \in \mathrm{Irr} \left(
\widetilde{\XXX \times \YYY}
\right)} 
\left( \deg_{\LLL \boxtimes \MMM} C \right) \delta_{\eta_{C}}
,
\end{align*}
where 
``$\mathrm{Irr}$''
stands for
the set of irreducible components
and $\eta_{A}$ denotes the point of the analytic space
corresponding to $A$
(cf. \S~\ref{bersp}).
Since $\widetilde{\XXX \times \YYY} = 
\tilde{\XXX} \times \tilde{\YYY}$,
we have naturally
$\mathrm{Irr} \left(
\widetilde{\XXX \times \YYY}
\right) 
=
\mathrm{Irr} 
\left(
\tilde{\XXX}
\right)
\times
\mathrm{Irr}
\left(
\tilde{\YYY}
\right)$.
If $C = A \times B$, then it is easy to see
\begin{align*}
\deg_{\LLL \boxtimes \MMM} C
=
\left(
\begin{matrix}
d + e
\\
d
\end{matrix}
\right)
(\deg_{\LLL} A)
\cdot
(\deg_{\MMM} B)
\end{align*}
and
$r_{\ast} \delta_{\eta_{C}} = \delta_{\eta_{A}} \times \delta_{\eta_{B}}$,
where $d := \dim X$ and $e := \dim Y$.
Accordingly, we have
\[
r_{\ast} 
\left(
\left(
\deg_{\LLL \boxtimes \MMM} C
\right)
\delta_{\eta_{C}} 
\right)
= 
\left(
\begin{matrix}
d + e
\\
d
\end{matrix}
\right)
\left(
\left(
\deg_{\LLL} A
\right)
\delta_{\eta_{A}} 
\right)
\times 
\left(
\left(
\deg_{\MMM} B
\right)
\delta_{\eta_{B}}
\right)
.
\]
Since
\begin{align*}
\deg_{L \boxtimes M} X \times Y
=
\left(
\begin{matrix}
d + e
\\
d
\end{matrix}
\right)
(\deg_{L} X)
\cdot
(\deg_{M} Y) 
,
\end{align*}
we thus obtain our formula in this case.

Now let us consider the general case.
Let $(\mathscr{X}_{n}, \mathscr{L}_{n})$
and
$(\mathscr{Y}_{n}, \mathscr{M}_{n})$ be
approximating sequences 
of models
of $\overline{L}$ and $\overline{M}$ respectively.
Then 
$\left( 
\XXX_{n} \times \YYY_{n}, 
\mathscr{L}_{n} \boxtimes \mathscr{M}_{n}
\right)$
is an approximating sequence of 
$\overline{L} \boxtimes \overline{M}$,
and
we have
\[
r_{\ast}
\left(
\mu_{\mathscr{L}_{n} \boxtimes \MMM_{n}}
\right)
=
\mu_{\mathscr{L}_{n}} \times \mu_{\MMM_{n}}
\]
as we have shown.
Taking the limit as $n \to + \infty$, we obtain our assertion.
\end{Pf}

\subsection{Non-degenerate strata} \label{nondeg}

We recall the notion of non-degenerate
strata here.
First of all, let us recall
the notion of stratification of a variety
(cf. \cite{berkovich1},
\cite{gubler3}).
Let $Z$ be a reduced scheme of finite type over a
field $k$.
Put $Z^{(0)} := Z$.
For $r \in \ZZ_{\geq 0}$, define
$Z^{(r+1)} \subset Z^{(r)}$ to be the complement of
the set of normal points of $Z^{(r)}$.
Then $Z^{(r+1)}$ is a proper closed subset of $Z^{(r)}$,
and we obtain a chain of closed subsets
\[
Z = Z^{(0)} \supsetneq Z^{(1)} \supsetneq
\cdots
\supsetneq Z^{(s-1)}
\supsetneq Z^{(s)}
=
\emptyset
\]
for some $s \in \NN$.
The irreducible component of $Z^{(r)} \setminus Z^{(r+1)}$
for any $r \in \ZZ_{\geq 0}$ is called a
\emph{stratum} of $Z$,
and the set of the strata of $Z$ is denoted by
$\str (Z)$.

We use here the same notations and conventions 
as those in \S~\ref{preliminary}.
Let $\XXX'$ be a strictly semistable formal scheme
(cf. \cite[5.1]{gubler3}).
Berkovich 
defined in \cite[\S 5]{berkovich4} the \emph{skeleton} $S ( \XXX ' )$.
It is a closed subset of $(\XXX')^{\an}$,
with
important properties as follows:
\begin{itemize}
\item
There is a continuous map
$\Val : (\XXX')^{\an} \to S ( \XXX' )$
which restricts to the identity on $S ( \XXX' )$.
\item
$S (\XXX')$ has a canonical structure
of metrized simplicial set;
there is
a family of metrized simplicial sets
$\{ \Delta_{S} \}_{S \in \str 
\left( \tilde{\XXX '} \right) }$
which
covers $S (\XXX')$.
\end{itemize}
Let $S$ be a stratum of $\tilde{\XXX'}$.
Let us describe $\Delta_{S}$ above a little more.
By the definition of strict semistability,
we can take a
$\pi \in \KK^{\circ \circ} \setminus \{ 0 \}$
and
an open subset $\UUU' \subset \XXX'$
with an \'etale morphism
\begin{equation*} 
\phi :
\UUU' \to
\Spf 
\KK^{\circ}
\langle 
x_{0}' , \ldots , x_{d}'
\rangle 
/
(x_{0}' \ldots x_{r}' - \pi)
\end{equation*}
such that $S \cap \tilde{\UUU'}$ dominates
$x_{0}' \ldots x_{r}' = 0$,
where $\KK^{\circ}
\langle 
x_{0}' , \ldots , x_{d}'
\rangle 
$
denote the Tate algebra.
Then 
we have an identification
\[
\left\{
\left. 
\mathbf{u'} \in \RR_{\geq 0}^{r+1}
\right|
u_{0}' + \cdots + u_{r}' =
v (\pi)
\right\}
\cong
\Delta_{S}
.
\]

Let $A$ be an abelian variety over $\KK$.
Recall that we have a continuous map 
$\overline{\val} : A^{\an} \to \RR^{n} / \Lambda$,
where $n = \dim A - b (A)$ and $\Lambda$ is a lattice
(cf. (\ref{tropicalization})).
Let $X \subset A$ be an irreducible closed subvariety of dimension $d$.
Let
$p : \EEE \to \AAA$ be a Mumford model of 
the Raynaud extension of $A$,
$\XXX \subset \AAA$ an admissible closed formal subscheme
with $\XXX^{\an} = X^{\an}$,
and let
$\XXX' \to \XXX$ be a semistable alteration,
that is,
$\XXX'$ is a strictly semistable formal scheme and
the morphism
$\XXX' \to \XXX$ is
a proper surjective generically finite morphism.
We denote by $f$ 
the composite
$
\XXX' \to \XXX \hookrightarrow \AAA$.
Gubler found in \cite[Proposition 5.11]{gubler3}
a unique continuous map
$\overline{f}_{\aff} : S (\XXX') \to \RR^{n} / \Lambda$
such that
$\overline{f}_{\mathrm{aff}} \circ \Val = \overline{\val} \circ f$.
Let $\YYY \subset \EEE$ be an open formal subscheme
such that $p (\YYY) = \XXX$.
We put $\YYY' := \XXX' \times_{\AAA} \YYY$,
and let $g$ be the composite
$
\YYY' \to \YYY \hookrightarrow \EEE
$.
We then
have a diagram
\begin{align*}
\begin{CD}
\YYY' @>{g}>> \EEE @>{q}>> \BBB\\
@V{p'}VV @VV{p}V \\
\XXX' @>{f}>> \AAA
,
\end{CD}
\end{align*}
in which the square commutes.
Let $S$ be a stratum of $\tilde{\XXX}'$.
Since $p'$ is surjective,
we can take an irreducible locally closed subset
$T \subset \tilde{\YYY'}$
such that $\tilde{p'} |_{T} : T \to S$
is dominant\footnote{It is automatically an open immersion.}.
With this notation,
we say
$S$
is
\emph{non-degenerate} with respect to $f$
if $\dim \overline{f}_{\mathrm{aff}} (\Delta_{S})
= \dim (\Delta_{S})$ and
$\dim
\left(
\tilde{q} \circ \tilde{g} (T)
\right)
=
\dim S$.
We also say $\Delta_{S}$ is \emph{non-degenerate}
with respect to $f$
if $S$ is non-degenerate with respect to $f$,
following Gubler's terminology
(cf. \cite[\S~6.3]{gubler3}).
The notion of non-degeneracy is well-defined
from
$S$ and $f$ ---
independent of any other choices.

Finally in this subsection,
we like to give a remark on the relationship 
between the dimension of abelian part and the dimension of the non-degenerate
strata.
To do that,
we assume in addition that the above $\YYY$ and hence $\YYY'$
are quasicompact.
Let $S$ be a non-degenerate stratum of $\tilde{\XXX}'$.
Then we have
\[
\max_{T' \in \str \left( \tilde{\YYY'} \right)}
\left\{
\dim 
\tilde{q} \left(
\tilde{g} \left( 
T' \right)
\right)
\right\}
\geq 
\dim S
\]
and hence
\[
b (X) = \dim \tilde{q} \left(
\tilde{g} \left( 
\tilde{\YYY'} \right)
\right)
=
\max_{T' \in \str \left( \tilde{\YYY '} \right)}
\dim 
\tilde{q} \left(
\tilde{g} \left( 
T' \right)
\right)
\geq 
\dim S
.
\]
Accordingly, we have
\addtocounter{Claim}{1}
\begin{align} \label{nondegeneratestr-b}
\dim \Delta_{S} 
= d - \dim S
\geq d - b (X)
.
\end{align}

\subsection{Minimum of the dimension of the
components of
the canonical measure}

Let $A$ be an abelian variety over $\overline{K}$
and let
$X \subset A$ be an irreducible closed subvariety of 
dimension $d$.
From now on,
we consider only canonical metrics
on line bundles
on abelian varieties,
hence
we let
$\overline{L}$
always
stand for
a line bundle $L$ with a
canonical  metric.

Let $v$ be a place of $\overline{K}$ 
and put $n := \dim A - b (A_{v})$. 
Since we have a continuous map
$\overline{\val} : A_{v} \to \RR^{n} / \Lambda$,
we can consider the tropicalization 
\[
\mu_{X_{v}, \overline{L}}^{\trop}
:=
\overline{\val}_{\ast} ( \mu_{X_{v}, \overline{L}})
\]
of the canonical measure,
which we call the \emph{tropical canonical measure}.
The measures
$\mu_{X_{v}, \overline{L}}$
and
$\mu_{X_{v}, \overline{L}}^{\trop}$
were studied in \cite{gubler3}.
We first recall the description obtained there:

\begin{citeTheorem}[The case of $L_{1} = \cdots = L_{d} = L$
in Theorem 1.1 in \cite{gubler3}]
\label{gublermeasure}
With the notation above,
suppose that $L$ is ample.
Then there are rational simplexes
$\overline{\Delta}_{1}, \ldots , \overline{\Delta}_{N}$
in $\RR^{n} / \Lambda$
with the following properties:
\begin{enumerate}
\renewcommand{\labelenumi}{(\alph{enumi})}
\item
$d - b (A_{v} ) \leq
\dim \overline{\Delta}_{j}
\leq
d$
for all $j = 1 , \ldots , N$.
\item
$X_{v}^{\trop} = \bigcup_{j = 1}^{N} \overline{\Delta}_{j}$.
\item
There are $r_{1}, \ldots , r_{N} > 0$
such that
\[
\mu_{X_{v}, \overline{L}}^{\trop} =
\sum_{j=1}^{N} r_{j} \delta_{\overline{\Delta}_{j}},
\]
where $\delta_{\overline{\Delta}_{j}}$ is the 
pushforward 
to $\RR^{n} / \Lambda$
of the 
canonical Lebesgue measure
on
$\overline{\Delta}_{j}$.
\end{enumerate}
\end{citeTheorem}

In general,
let
$\mu$
be a measure on a polytopal subset of $\RR^{n} / \Lambda$
of form
\[
\mu = 
\sum_{i = 1}^{N} r_{i} \delta_{\overline{\Delta}_{i}}
\quad
(r_{i} > 0)
.
\]
Then we define $\sigma ( \mu )$
by
\begin{align*} 
\sigma (\mu)
:=
\min_{i}
\{
\dim \overline{\Delta}_{i}
\}
.
\end{align*}

Let $\XXX$ be the closure of $X_{v}$ in a
Mumford model of $A_{v}$.
We can take a semistable alteration $f : \XXX ' \to 
\XXX$ of a model $\XXX$ of $X_{v}$
by virtue of \cite[Theorem 6.5]{dejong}.
We can write
\begin{align*}
c_{1} (f^{\ast} \overline{L})^{\wedge d} =
\sum_{S}
r_{S} \delta_{\Delta_{S}}
\end{align*}
by \cite[Corollary 6.9]{gubler3},
where $S$ ranges over all the non-degenerate strata 
of $\tilde{\XXX'}$ with respect
to $f$, and $r_{S}$ is positive.
By \cite[Propositions 3.9 and 5.11]{gubler3},
we have
\[
\overline{\mathrm{val}}_{\ast}
\left(
c_{1} ( \overline{L})^{\wedge d}
\right)
=
(\deg f)
(\overline{f}_{\mathrm{aff}})_{\ast}
(c_{1} (f^{\ast} \overline{L})^{\wedge d})
.
\]
Therefore we can write
\begin{align*}
\mu_{X_{v}, \overline{L}}^{\trop} 
=
\sum_{S}
r'_{S} \delta_{\overline{f}_{\mathrm{aff}}(\Delta_{S})}
\end{align*}
for some $r'_{S} > 0$.
Since $\Delta_{S}$ is non-degenerate,
we have
$\dim \Delta_{S} = \dim f_{\mathrm{aff}}(\Delta_{S})$.
We see therefore
\addtocounter{Claim}{1}
\begin{equation} 
\label{tropcanonical}
\begin{split}
\sigma 
\left( 
\mu_{X_{v}, \overline{L}}^{\trop}
\right)
&=
\min
\left\{
\dim \Delta_{S}
\left|
\text{$S \in \str 
\left(
\tilde{\XXX'}
\right)$ is 
non-degenerate with respect to $f$}
\right.
\right\}
\\
&=
d -
\max
\left\{
\dim S
\left|
\text{$S \in \str 
\left(
\tilde{\XXX'}
\right)$ is 
non-degenerate with respect to $f$}
\right.
\right\}
,
\end{split}
\end{equation}
and combining it with
(\ref{nondegeneratestr-b}),
we obtain
\addtocounter{Claim}{1}
\begin{align} \label{tropcanonical-b}
\sigma 
\left( 
\mu_{X_{v}, \overline{L}}^{\trop}
\right)
\geq
d - b (X_{v})
.
\end{align}



\section{Conclusions}

\subsection{Special subvariety and the dimension of
the abelian part} \label{abpt}

We begin this section with
the following assertion.

\begin{Proposition} \label{compprop}
Let $X$ be a special subvariety of $A$.
Then we have the following.
\begin{enumerate}
\item
$\dim X / G_{X} = b ((X / G_{X})_{v})$
for any place $v$.
\item
Suppose
$\dim X / G_{X} \geq b ((A / G_{X})_{v})$.
Then
there is a special point $\sigma$ with
$X = G_{X} + \sigma$,
that is, $X$ is an abelian subvariety up to a special point.
In particular,
we have
$\dim X / G_{X} = b ((X / G_{X})_{v}) = 0$.
\end{enumerate}
\end{Proposition}

\begin{Pf}
Taking the quotient by $G_{X}$, we may assume $G_{X} = 0$
(cf. Proposition~\ref{specialsubvariety}),
and further,
taking the translate of $X$ by a torsion point if necessary,
we may assume
that there is a closed subvariety 
$Y' \subset A^{\overline{K}/k}$
such that 
\[
\Tr_{A}^{\overline{K}/k} ({Y'}_{\overline{K}})
= X.
\]

We put $\KK := \overline{K}_{v}$.
By the existence of the N\'eron model
and the semistable reduction theorem,
we can take
a semi-abelian scheme $\mathcal{A}$
over $\KK^{\circ}$
and a homomorphism
\[
\tau :
A^{\overline{K}/k} \times_{\Spec k}
\Spec \KK^{\circ}
\to
\mathcal{A}
\]
extending 
$\Tr_{A}^{\overline{K}/k}$.
Since $A^{\overline{K}/k} \times_{\Spec k}
\Spec \KK^{\circ}$ is proper over
$\Spec \KK^{\circ}$, we see
that $\tau$ is proper.
Let $\tilde{\tau}$ denote the 
restriction
of $\tau$ to the special fiber.

\begin{Claim} \label{taufinite}
$\tilde{\tau}$ is a finite morphism.
\end{Claim}

\begin{Pf}
Let $\mathcal{Z}$ be the scheme theoretic image of $\tau$ in $\mathcal{A}$.
Since $\tau$ is proper,
the morphism $A^{\overline{K}/k} \times_{\Spec k}
\Spec \KK^{\circ}
\to \mathcal{Z}$ is proper and surjective.
Since $\tau$ is finite over 
$\mathcal{Z}_{\KK} := \mathcal{Z} \times_{\Spec \KK^{\circ}} \Spec \KK$
by Lemma~\ref{finiteinsep},
we have
$\dim \mathcal{Z}_{\KK} = \dim A^{\overline{K}/k}$.
Therefore
we find
$\dim \tilde{\mathcal{Z}} = \dim A^{\overline{K}/k}$
by Chevalley's theorem (\cite[Th\'eor\`eme 13.1.3]{ega28}),
from
which we see that
$\tilde{\tau} : 
A^{\overline{K}/k} \to \tilde{\mathcal{Z}}$ is 
a 
generically finite proper surjective morphism.
Since $\tilde{\tau}$ is a surjective homomorphism of group schemes,
its fiber is equidimensional.
Thus we conclude that
$\tilde{\tau}$ is finite.
\end{Pf}

\begin{Claim}
We have
$b (X_{v}) 
\geq  
\dim Y'$
and
$b (A_{v}) \geq \dim A^{\overline{K}/k}$.
\end{Claim}

\begin{Pf}
We here fix a Mumford model $p : \EEE \to \AAA$
of the uniformization $E \to A_{v}$.
Recall that,
for our $A_{v}$,
we have a unique exact sequence
\begin{align*}
\begin{CD}
1 @>>> {\mathscr{T}}^{\circ}
@>>> {\EEE}^{\circ} @>{{q}^{\circ}}>>
{\BBB} @>>> 0
\end{CD}
\end{align*}
as
(\ref{formalraynaud2}),
and an isomorphism $p^{\circ} : \EEE^{\circ} \to \AAA^{\circ}$
as in (\ref{modelofopensubgroups}).
According to 
the construction of Raynaud extension in
\cite[\S 1]{BL},
we can identify
$\AAA^{\circ}$ with
the formal completion $\hat{\mathcal{A}}$ of $\mathcal{A}$.
Via this identification,
we regard
the reduction $\tilde{\tau}$
as
a 
homomorphism from $A^{\overline{K}/k}$ to $\tilde{\AAA}^{\circ}$.
Since $
(\tilde{p}^{\circ})^{-1}
\left(
\tilde{\tau} \left( A^{\overline{K}/k} \right)
\right)
$
is 
proper over $k$
and $\Ker \tilde{q}^{\circ}$ is affine,
we see that
$\tilde{q}^{\circ} |_{(\tilde{p}^{\circ})^{-1}
\left(
\tilde{\tau} \left( A^{\overline{K}/k} \right)
\right)}$
is finite.
Therefore
we find $\tilde{q}^{\circ} \circ (\tilde{p}^{\circ})^{-1} \circ \tilde{\tau}$
is also finite by Claim~\ref{taufinite}.

Now the second inequality
of our claim is obvious;
\[
b (A_{v}) = \dim \tilde{\BBB} \geq
\dim
\tilde{q}^{\circ} 
\left(
(\tilde{p}^{\circ})^{-1}
\left( \tilde{\tau}  \left( A^{\overline{K}/k} \right) \right)
\right)
=
\dim A^{\overline{K}/k}
.
\]
Let us show the first inequality.
Let $\XXX$ be the closure of $X$ in $\AAA$,
and put
$\XXX^{\circ} := \XXX \cap \AAA^{\circ}$.
We put 
\[
\YYY^{\circ} := (p^{\circ})^{-1} (\XXX^{\circ})
\subset \EEE^{\circ}
,
\]
and let $\YYY \subset \EEE$ be a quasicompact open
formal subscheme
such that $\YYY^{\circ} \subset \YYY$ and
$p (\tilde{\YYY}) = \XXX$.
We can see that
the special fiber of the closure of $X$ in $\mathcal{A}$ coincides with
$\tilde{\XXX}^{\circ}$
via the identification $\tilde{\mathcal{A}} = \tilde{\AAA}^{\circ}$.
Taking account of $ \tilde{\tau} (Y') \subset
\tilde{\XXX}^{\circ}$,
which comes from  our assumption,
we have
\begin{multline*}
b(X_{v}) = \dim \tilde{q} 
\left(
\tilde{\YYY}
\right)
\geq 
\dim
\tilde{q}^{\circ} 
\left( \tilde{\YYY}^{\circ} 
\right)
=
\dim
\tilde{q}^{\circ} \left(
\left(
\tilde{p}^{\circ}
\right)^{-1}
\left( \tilde{\XXX}^{\circ} 
\right)
\right)
\\
\geq
\dim
\tilde{q}^{\circ} 
\left(
\left(
\tilde{p}^{\circ}
\right)^{-1}
\left( \tilde{\tau} 
\left( Y'
\right)
\right)
\right)
=
\dim (Y')
\end{multline*}
as required.
\end{Pf}

Accordingly,
we have
\[
\dim Y' = \dim X
\geq b (X_{v}) \geq \dim Y'
\]
and hence $\dim X = b (X_{v})$, which proves (1).
To show (2), we suppose 
further $\dim X  \geq b (A_{v})$.
Then we have
\[
\dim Y' = \dim X
\geq b (A_{v}) \geq \dim A^{\overline{K} / k}
\geq \dim Y'
,
\]
which says
$Y' = A^{\overline{K} / k}$
and 
$X$ is an abelian subvariety.
Since
we have
$G_{X} = 0$ by our assumption,
we find $X = 0$
as required.
\end{Pf}

\begin{Remark} \label{spcial-abeliansub}
Suppose that $X$ is a special subvariety of $A$ over $\overline{K}$ and 
that there is a place $v \in M_{\overline{K}}$
at which $A$ is totally degenerate.
Then it immediately follows from
Theorem~\ref{thmofgubler} and
Corollary~\ref{special-small}
that
$X$ is a torsion subvariety,
but we can show this fact directly (without using Gubler's theorem).
In fact,
we have 
$b (( A / G_{X}))_{v} = 0$ by Lemma~\ref{quotientb}
and hence
$X$ is the translate of $G_{X}$ by a special point
by Proposition~\ref{compprop} (2).
Any special point is a torsion point since $A_{v}$ is totally degenerate,
and hence
we conclude that $X$ is
a torsion subvariety in this case.
We can also show by a similar argument that
a special subvariety is
an abelian subvariety up to translation by a special point
in the case where there exists a place $v$ with $b (A_{v}) = 1$.
\end{Remark}

\subsection{First main result}

According to Proposition~\ref{compprop} (1),
an irreducible closed subvariety with
$\dim (X / G_{X}) > b ((X / G_{X})_{v})$
for some $v$ is not a special subvariety.
If the geometric Bogomolov conjecture holds true,
then such a closed subvariety should not
have dense small points.
In fact, it is our main assertion:

\begin{Theorem} [cf. Theorem~\ref{thoremint}]
\label{maintheorem}
Let $A$ be an abelian variety over $\overline{K}$
and
let $X$ be an irreducible closed subvariety of $A$.
Let $G_{X} \subset A$ be the stabilizer of $X$.
Suppose 
$\dim (X / G_{X}) > b ((X / G_{X})_{v})$
for some place $v$.
Then $X$ does not have
dense small points.
\end{Theorem}

\begin{Pf}
We argue by contradiction.
Suppose that there exists a counterexample
$X$ to Theorem~\ref{maintheorem},
that is,
$X$
has dense small points but is not a special subvariety.
Then
the closed subvariety $X / G_{X} \subset A / G_{X}$
is again a counterexample
by virtue of Lemma~\ref{image-counterexample}
and Proposition~\ref{specialsubvariety}.
Accordingly
we may assume
$G_{X} = 0$
and 
our assumption 
in the theorem says
$d := \dim X > b (X_{v})$.
Since $G_{X} = 0$,
there exists an integer $N > 0$ such that
\[
\alpha : X^{N} \to A^{N-1}
,
\quad
(x_{1}, \ldots , x_{N})
\mapsto
(x_{2} - x_{1}, \ldots , x_{N} - x_{N-1})
\]
is generically finite
(cf. \cite[Lemma~3.1]{zhang2}).
We put $X' := X^{N}$
and $Y := \alpha (X')$.
The closed subvariety
$X' \subset A^{N}$
also has dense small points
by Lemma~\ref{product-dense}. 

Let $L$ and $M$ be even ample line bundles on $X$
and $Y$ respectively.
Then the line bundle $L' := L^{\boxtimes N}$
on $A^{N}$ is even and ample.
Let $\mu$ and $\nu$ be the 
tropical canonical
measures on 
$(X'_{v})^{\trop} = (X_{v}^{\trop})^{N}$ 
and
$Y_{v}^{\trop}$
arising from
$L'$ and $M$ respectively.
We simply write $\hat{h}_{X'}$ and 
$\hat{h}_{Y}$ for the canonical heights
associated with $L'$ and $M$ respectively.
Since $X'$ has dense small points,
we can find 
a generic net $(P_{m})_{m \in I}$,
where $I$ is a directed set,
such that
$\lim_{m} \hat{h}_{X'} (P_{m}) = 0$.
We call such an net a 
\emph{generic net of small points}.
The image $(\alpha (P_{m}))_{m \in I}$
is also a generic net 
of small points
of $Y$.
Then by using the equidistribution theorem
\cite[Theorem 1.2]{gubler4},
we can obtain
\[
(\overline{\alpha}^{\trop})_{\ast} \mu = \nu
\]
in the usual way
(cf. \cite[Proof of Theorem 1.1]{gubler2}),
where $\overline{\alpha}^{\trop} : 
(X'_{v})^{\trop} \to Y_{v}^{\trop}$
is the
map between tropical varieties associated to
$\alpha$.
(In Gubler's article, it is denoted by 
$\overline{\alpha}_{{\val}}$.)

Let us take Mumford models
$\AAA_{1}$ of $A_{v}^{N}$ and $\AAA_{2}$ of $A_{v}^{N-1}$
such that $\alpha : X_{v}' \to Y_{v}$
extends to the morphism of models
$h : \XXX' \to \YYY$,
where $\XXX'$ is the closure of $X_{v}'$
in $\AAA_{1}$
and $\YYY$ is that of $Y$ in $\AAA_{2}$.
Let $f : \XXX'' \to \XXX'$
be a strictly semistable alteration.
Then $g := h \circ f$
is also a strictly semistable alteration for $\YYY$
since $h$ is a generically finite surjective morphism.
Let $S$ be a stratum of $\tilde{\XXX''}$.
Then,
we immediately see from the definition of
non-degeneracy
that
$S$ is non-degenerate with respect to $f$
if so is $S$ with respect to $g$.
In particular we have
\begin{multline*}
\max 
\{
\dim S
\mid
\text{$S$ is non-degenerate with respect to $f$}
\}
\\
\geq
\max 
\{
\dim S
\mid
\text{$S$ is non-degenerate with respect to $g$}
\}
,
\end{multline*}
and 
we find
\addtocounter{Claim}{1}
\begin{align} \label{keyinequality}
\sigma (\mu ) 
\leq
\sigma (\nu)
\end{align}
by (\ref{tropcanonical}).%
\footnote{We have $\sigma (\mu ) 
= 
\sigma (\nu)$ in fact.}

Let us write
\[
\mu_{X_{v}, \overline{L}}^{\trop}
=
\sum_{j = 1}^{N} r_{j} \delta_{\overline{\Delta}_{j}}
,
\quad
(r_{j} > 0)
\]
as in Theorem~\ref{gublermeasure}.
Renumbering them if necessary,
we may assume
$\dim \overline{\Delta}_{1} = 
\sigma 
\left(
\mu_{X_{v}, \overline{L}}^{\trop}
\right)$.
Since $d > b (X_{v})$ by our assumption,
we have
$\dim \overline{\Delta}_{1} > 0$ 
by (\ref{tropcanonical-b}).
Taking account of Lemma~\ref{product-trop}
and Proposition~\ref{productmeasure},
we can write
\[
\mu
=
\sum_{j_{1}, \ldots , j_{N}}
r_{j_{1}} \ldots r_{j_{N}}
\left(
\delta_{\overline{\Delta}_{j_{1}}} \times \cdots
\times \delta_{\overline{\Delta}_{j_{N}}}
\right)
=
\sum_{j_{1}, \ldots , j_{N}}
r_{j_{1}} \ldots r_{j_{N}}
\left(
\delta_{\overline{\Delta}_{j_{1}} \times \cdots
\times \overline{\Delta}_{j_{N}}}
\right)
.
\]
The
coefficients
in the summation 
are all positive,
and we have
\[
\dim {\overline{\Delta}_{1}}^{N} = 
\sigma (\mu)
=
N \sigma
\left(
\mu_{X_{v}, \overline{L}}^{\trop}
\right) > 0
.
\]

Since $\alpha$ contracts
the diagonal of $X'$ to the origin of 
$A^{N-1}$,
we see
$\overline{\alpha}^{\trop}$ also contracts
that of
${\overline{\Delta}_{1}}^{N}$
to $\overline{\mathbf{0}}$.
Therefore,
there exists a $\sigma (\mu)$-dimensional simplex 
$\overline{\Delta} \subset {\overline{\Delta}_{1}}^{N}$ 
such that
$\dim {\overline{\alpha}}^{\trop} (\overline{\Delta})
< \sigma (\mu)$.
On the other hand,
we have $\nu (\overline{\tau}) = 0$
for any simplex $\tau$ of dimension less than
$\sigma (\mu)$
by (\ref{keyinequality}),
which says
$\nu (\overline{\alpha}^{\trop} (\overline{\Delta}) ) = 0$
in particular.
On the other hand,
since $(\overline{\alpha}^{\trop})_{\ast} \mu = \nu$,
we have
\[
\nu (\overline{\alpha}^{\trop} (\overline{\Delta}) )
=
\mu
\left(
(\overline{\alpha}^{\trop})^{-1}
\left(
\overline{\alpha}^{\trop} (\overline{\Delta})
\right)
\right)
\geq
\mu ( \overline{\Delta})
> 0
.
\]
That is a contradiction,
and thus
we complete the proof.
\end{Pf}

\subsection{Further results} \label{furtherresults}

In this subsection, we use a result of the appendix
by Gubler
to show some results concerning the geometric Bogomolov conjecture,
though we have not used the appendix so far.

\begin{Theorem} \label{addthm}
Let $X$ be an irreducible subvariety of an abelian variety $A$.
Suppose there exists a place $v$ such that 
$\dim X / G_{X} \geq b ((A / G_{X})_{v})$.
Then $X$ does not have dense small points
unless it is a special subvariety.
\end{Theorem}

\begin{Pf}
We argue by contradiction.
Suppose that $X$ has dense small points.
As usual,
we may assume $G_{X} = 0$ by taking the quotient,
and
$\dim X > 0$ by Lemma~\ref{zerodimensional}.
If $\dim X > b(X_{v})$, we are done by Theorem~\ref{maintheorem}.
Consider the case $\dim X \leq b(X_{v})$.
Then $\dim X = b(X_{v})$ and hence $b(X_{v}) \geq b (A_{v})$ by our assumption,
which concludes $b(X_{v}) = b (A_{v})$.
Note that $b (A_{v}) > 0$ in this situation.

Consider the morphism $\alpha : X^{N} \to A^{N-1}$
as in the proof of Theorem~\ref{maintheorem},
and put $X' := X^{N}$ and $Y := \alpha (X')$
as before.
Recall that $\alpha : X' \to Y$ is 
a generically finite surjective morphism.
We have
\[
b (X'_{v}) = N b(X_{v}) = N b (A_{v}) > (N-1) b (A_{v}) \geq b (Y_{v})
.
\]

Let $\mu$ and $\nu$ be the
canonical measures on $X'_{v}$
and $Y_{v}$ respectively,
which are the same ones as in the proof of Theorem~\ref{maintheorem}.
Then by Gubler's result Corollary~\ref{tropical canonical measure},
we obtain
\addtocounter{Claim}{1}
\begin{align} \label{inequality1}
\sigma 
\left( \mu \right)
=
\dim X' - b (X'_{v})
<
\dim Y - b (Y_{v}) =
\sigma
\left(
\nu
\right)
.
\end{align}
On the other hand, $X'$ also has dense small points.
Therefore we have
a generic net of small points,
and
the image of this generic net by $\alpha$ is also
a generic net of small points of $Y$.
By the equidistribution theorem \cite[Theorem 1.1]{gubler4},
we then have
\[
\left(
\overline{\alpha}^{\trop}
\right)_{\ast}
\mu = \nu
,
\]
which implies
\[
\sigma (\nu) = \sigma 
\left(
\left(
\overline{\alpha}^{\trop}
\right)_{\ast}
\mu
\right)
\leq \sigma ( \mu )
.
\]
That however contradicts (\ref{inequality1}).
\end{Pf}

\begin{Corollary} \label{b(A/G)leq1}
Let $X$ and $A$ be as above.
Suppose $b ((A / G_{X})_{v}) \leq 1$.
If $X$ has dense small points, then it is a special subvariety.
\end{Corollary}

\begin{Pf}
Suppose that
$X$ has
dense small points.
Since $b( (A/G_{X})_{v}) \leq 1$,
we have $\dim X / G_{X} = 0$ by Theorem~\ref{addthm}.
By Lemma~\ref{zerodimensional},
we then conclude that $X$ is a special subvariety.\footnote{As
we saw in Remark~\ref{spcial-abeliansub},
$X$ 
is the translate of $G_{X}$
by a special point of $A$ in fact.}
\end{Pf}

As a further corollary, we can 
obtain
the following assertion.
In the case of $b (A_{v}) = 0$, it is \cite[Theorem 1.1]{gubler2}.

\begin{Corollary}[cf. Theorem~\ref{GBC:bleq1int}] \label{GBC:bleq1}
Let $A$ be an abelian variety.
Suppose that there exists a place 
$v$ such that 
$b (A_{v}) \leq 1$.
Then the geometric Bogomolov conjecture holds for $A$.
\end{Corollary}

\begin{Pf}
This assertion immediately follows from 
Corollary~\ref{b(A/G)leq1}
since we have
$b( (A/G_{X})_{v}) \leq 1$ by Lemma~\ref{quotientb}.
\end{Pf}



\renewcommand{\thesection}{Appendix by Walter Gubler} 
\renewcommand{\theTheorem}{A.\arabic{Theorem}}
\renewcommand{\theClaim}{A.\arabic{Theorem}.\arabic{Claim}}
\renewcommand{\theequation}{A.\arabic{Theorem}.\arabic{Claim}}
\renewcommand{\theProposition}
{A.\arabic{Theorem}.\arabic{Proposition}}
\renewcommand{\theLemma}{A.\arabic{Theorem}.\arabic{Lemma}}
\setcounter{section}{0}
\renewcommand{\thesubsection}{A.\arabic{subsection}}
\renewcommand{\theart}
{A.\arabic{Theorem}.\arabic{}}



\section{The minimal dimension of a canonical measure}

In Yamaki's proof of Theorem~\ref{thoremint}, the main point was to deduce
the lower bound 
$d - b(X_{v})$
for the minimal dimension of the tropical
canonical measure
(cf. (\ref{tropcanonical-b})). 
We will show in this appendix that the minimal
dimension is in fact equal to
$d - b(X_{v})$
and that this holds also for a
canonical measure on X. This is interesting as such measures play an
important role in non-archimedean analysis.

Let $K$ be a field with a discrete valuation $v$ and let $\kdop=\cdop_K$ be a minimal algebraically closed field which is complete with respect to a valuation extending $v$. The valuation ring of $\kdop$ is denoted by $\kcirc$. 

We consider an irreducible $d$-dimensional closed subvariety $X$ of an abelian variety $A$ defined over $\overline{K}$. We will recall in \ref{canonical subset} that the Berkovich analytic space $\Xan$ over $\kdop$ associated to $X$ has a canonical piecewise linear subspace $S_{X}$ which is the support of every canonical measure on $X$. Let $b(X)$ be the dimension of the abelian part of $X$ (see \S~\ref{abelianpart}). We will also use the uniformization $p:E \rightarrow \Aan = E/M$ from the Raynaud extension and the corresponding tropicalization maps $\val:E \rightarrow \rdop^n$ and $\valbar: \Aan \rightarrow \rtor$ (see \S~\ref{RtM}).

The goal of this appendix is to show the following result.
\begin{Theorem} \label{canonical measure}
There are rational simplices $\Delta_1, \dots , \Delta_{N}$ in $S_X$ with the following five properties:
\begin{itemize}
\item[(a)] For $j=1, \dots,N$, we have $\dim(\Delta_j) \leq d$.
\item[(b)] $S_X= \bigcup_{j=1}^N \Delta_j$.
\item[(c)] The restriction of $\valbar$ to $\Delta_j$ induces a linear isomorphism onto a simplex $\Deltabar_j$ of $\rtor$. 
\item[(d)] For canonically metrized line bundles $\overline{L_1}, \dots, \overline{L_d}$ on $A$, there are $r_j \in \rdop$ with
$$c_1(\overline{L_1}|_X) \wedge \dots \wedge c_1(\overline{L_d}|_X)
= \sum_{j=1}^{N} r_j \cdot \delta_{\Delta_j},$$
where $\delta_{\Delta_j}$ is the pushforward of the Lebesgue measure on $\Delta_j$ normalized by $\delta_{\Delta_j}(\Delta_j)=(\dim(\Delta_j)!)^{-1}$.
\item[(e)] If all line bundles in (d) are ample, then $r_j>0$ for all $j \in \{1, \dots , N\}$.
\end{itemize}
For any such covering of $S_X$, we have $\min\{\dim( \Delta_j)\mid j=1, \dots, N\} = d-b(X)$.
\end{Theorem}

The proof will be given in \ref{proof of canonical measure}.

\begin{Corollary} \label{tropical canonical measure}
Let $\Deltabar_1, \dots, \Deltabar_N$ be the components of the tropical canonical measure $\mu^{\rm trop}_{\Xan, \overline{L}}$ considered in Theorem~\ref{gublermeasure}. Then we have 
$$\min_{j=1,\dots, N}\dim( \Deltabar_j) = d-b(X).$$
\end{Corollary}

\begin{Pf} The tropical canonical measure satisfies
 $$\mu^{\rm trop}_{\Xan, \overline{L}}= 
 \frac{1}{\deg_{L} X} \valbar_*\left( c_1(\overline{L}|_X)^{\wedge d}  \right)$$
and hence Corollary \ref{tropical canonical measure} follows from Theorem \ref{canonical measure}. 
\end{Pf}

\begin{art} \label{setup} \rm 
Let $\Acal_0$ be the Mumford model of $A$ over $\kcirc$ associated to a rational polytopal decomposition $\overline{\Ccal_0}$  of $\rtor$. We denote the closure of $\Xan$ in $\Acal_0$ by $\Xcal_0$ which is a formal $\kcirc$-model of $\Xan$. It follows 
from de Jong's alteration theorem that there is a proper surjective morphism $\varphi_0:\Xcal' \rightarrow \Xcal_0$ from a strictly semistable formal scheme $\Xcal'$ over $\kcirc$ whose generic fibre  is an irreducible $d$-dimensional proper algebraic variety $X'$ (see \cite[6.2]{gubler3}). The generic fibre of $\varphi_0$ is denoted by $f$.
\end{art}

\begin{art} \label{canonical subset} \rm 
The {\it canonical subset} $S_X$ of $\Xan$ is defined as the support of a canonical measure $c_1(\overline{L_1}|_X) \wedge \dots \wedge c_1(\overline{L_d}|_X)$. Similarly as in \cite[Remark 6.11]{gubler3}, the definition of $S_X$ does not depend on the choice of the canonically metrized ample line bundles  $\overline{L_1}, \dots, \overline{L_d}$   of $A$. By \cite[Theorem 6.12]{gubler3}  $S_X$ is a rational piecewise linear space. The piecewise linear structure is characterized by the fact that the restriction of $f$ to the union of all canonical simplices which are non-degenerate with respect to $f$ induces a piecewise linear map onto $S_X$ with finite fibres. This structure does not depend on the choice of $\Acal_0$ and $f$ in \eqref{setup}. 
\end{art}

\begin{Theorem} \label{semistable canonical measure}
Let $\varphi:\Xcal' \rightarrow \Xcal_0$ be a strictly semistable alteration as in \ref{setup} with generic fibre $f:(X')^{\rm an} \rightarrow \Xan$. Then there is a $b(X)$-dimensional stratum $S$ of $\tilde{\Xcal}'$ such that the canonical simplex $\Delta_S$ of $S(\Xcal')$ is non-degenerate with respect to $f$.
\end{Theorem}

\begin{Pf}
We use the same method as in the proofs of  Theorem 6.7 and Lemma 7.1 in \cite{gubler3}.  Let $\Sigmabar$ be the collection of simplices of $X^{\rm trop}=\valbar(\Xan)$  given by  $\overline{f}_{\rm aff}(\Delta_S)$ together with all their closed faces where $S$ ranges over all strata of $\tilde{\Xcal}'$. There is a rational polytopal decomposition $\overline{\Ccal_1}$  of $\rtor$  which is transversal to $\Sigmabar$, i.e. $\Deltabar \cap \overline{\sigma}$ is either empty or of dimension  $\dim(\Deltabar)+\dim(\overline{\sigma}) - n$ for all $\Deltabar \in \overline{\Ccal_1}$ and $\overline{\sigma} \in \Sigmabar$. Note that the existence of such a transversal $\overline{\Ccal_1}$ is much easier than the construction in \cite[Lemma 6.5]{gubler3}, and no extension of the base field is needed here.

We consider the polytopal decomposition  $\Ccalbar:=\{\Deltabar_0 \cap \Deltabar_1 \mid  \Deltabar_0 \in \overline{\Ccal}_0, \Deltabar_1 \in \overline{\Ccal}_1\}$ which is the coarsest refinement of $\overline{\Ccal_0}$ and  $\overline{\Ccal_1}$. Let $\Acal_1, \Acal$ be the Mumford models associated to $\overline{\Ccal_1}$ and $\overline{\Ccal}$. Then we get  the following commutative diagram of canonical morphisms of  formal schemes over $\kcirc$:
\begin{equation*} 
\begin{CD} 
\Xcal'' @>\varphi>> \Acal @>{\iota_1}>> \Acal_1\\
@VV{\iota'}V    @VV{\iota_0}V\\
\Xcal' @>{\varphi_0}>> \Acal_0
\end{CD}
\end{equation*}
Here  the formal scheme $\Xcal''$ with reduced  special fibre is determined by the fact that the rectangle is cartesian on the level of formal analytic varieties (see \cite[5.17]{gubler3}). 

Let $\Ecal_0,\Ecal_1,\Ecal$ be the $\kcirc$-models of the uniformization $E$ associated to $\Ccal_0,\Ccal_1,\Ccal$ (see \S~\ref{RtM2}). For $i=1,2$, let $\iota_i':\Ecal \rightarrow \Ecal_i$ be the unique morphism extending the identity on the generic fibre. By construction, we have $\Acal_i:=\Ecal_i/M$  and $\Acal=\Ecal/M$ with quotient morphisms $p_i$ and $p$. 
The homomorphism $q:E \rightarrow B$ from the Raynaud extension is the generic fibre of  unique morphisms $q_i:\Ecal \rightarrow \Bcal$ and $q:\Ecal \rightarrow \Bcal$. Let  $\Xcal_1$ (resp. $\Xcal$) be the closure of $X$ in $\Acal_1$ (resp. $\Acal$) and let  $\Ycal_1:=p_1^{-1}(\Xcal_1)$, $\Ycal:=p^{-1}(\Xcal)$. By definition of $b(X)$, there is an irreducible component $W_1$ of $\tilde{\Ycal}_1$ with
\addtocounter{Claim}{1}
\begin{equation} \label{equation 1}
 \dim \tilde{q}_1(W_1) = b(X).
\end{equation} 
Since $\Ycal_1=\iota_1'(\Ycal)$, there is an irreducible component $W$ of $\tilde{\Ycal}$ with $W_1=\tilde{\iota}_1'(W)$. By \cite[Propositions 5.7 and 5.13]{gubler3}, there is a bijective correspondence between the vertices of the polytopal subdivision
$$\Dcal:=\{\Delta_S \cap \overline{f}_{\rm aff}^{-1}(\Deltabar) \mid \text{$S$ stratum of $\tilde{\Xcal}'$ and $\Deltabar  \in \Ccalbar$}\}$$
of the skeleton $S(\Xcal')$ and the $d$-dimensional strata of $\tilde{\Xcal}''$. Since $\tilde{p}$ is a local isomorphism, it is clear that $\tilde{p}(W)$ is an irreducible component of $\tilde{\Xcal}$. Using the fact that $\tilde{\varphi}$ is a proper surjective morphism onto $\tilde{\Xcal}$, there is a $d$-dimensional stratum $R$ of $\tilde{\Xcal}''$ with $\tilde{\varphi}(R)$ dense in $\tilde{p}(W)$.  Let $u'$ be the vertex of $\Dcal$ corresponding to $R$ and let $S$ be the unique stratum of $\tilde{\Xcal}'$  with $u'$ contained in the relative interior $\relint(\Delta_S)$.

By \cite[Lemma 5.15]{gubler3}, we have $\tilde{\iota}'(R)=S$, the map $\tilde{\varphi}_0:S \rightarrow \tilde{\Acal}_0=\tilde{\Ecal}_0/M$ has a lift $\tilde{\Phi}_0:S \rightarrow \tilde{\Ecal}_0$ and there is a unique lift $\tilde{\Phi}: R \rightarrow \tilde{\Ecal}$ of $\tilde{\varphi}:R \rightarrow \tilde{\Acal}=\tilde{\Ecal}/M$ with $\tilde{\Phi}_0 \circ \tilde{\iota}' = \tilde{\iota}_0' \circ \tilde{\Phi}$ on $R$. The lift $\tilde{\Phi}_0$ is unique up to $M$-translation and hence we may fix it by requiring that $\tilde{\Phi}(R)$ is dense in $W$. It follows that 
\addtocounter{Claim}{1}
\begin{equation} \label{equation 2}
 \tilde{q}_0(\tilde{\Phi}_0(S))=  \tilde{q}_0 \circ\tilde{\Phi}_0 \circ\tilde{\iota}'(R)=  \tilde{q}_0 \circ\tilde{\iota}_0' \circ \tilde{\Phi} (R)= 
\tilde{q}_1 \circ\tilde{\iota}_1' \circ \tilde{\Phi} (R)
\end{equation}
is dense in $\tilde{q}_1(W_1)$. By \eqref{equation 1}, we get 
\addtocounter{Claim}{1}
\begin{equation} \label{equation 3}
\dim \tilde{q}_0(\tilde{\Phi}_0(S)) = b(X).
 \end{equation}
Since $u'$ is a vertex of $\Dcal$ contained in $\relint(\Delta_S)$, it is clear that
\addtocounter{Claim}{1}
\begin{equation} \label{equation 4}
\dim \overline{f}_{\rm aff}(\Delta_S) = \dim \Delta_S
\end{equation}
(see also the argument after (25) in \cite[Remark 5.17]{gubler3}). There is a unique $\Deltabar_1 \in \overline{\Ccal_1}$ with  $\overline{f}_{\rm aff}(u') \in \relint(\Deltabar_1)$. Since $\overline{f}_{\rm aff}(u')$ is also contained in $\overline{f}_{\rm aff}(\Delta_S) \in \Sigmabar$, the transversality of  $\overline{\Ccal_1}$ and $\Sigmabar$ yields
\addtocounter{Claim}{1}
\begin{equation} \label{equation 5}
\codim \Deltabar_1  \leq \dim \overline{f}_{\rm aff}(\Delta_S) = \dim \Delta_S.
\end{equation}
By \cite[Proposition 5.14]{gubler3},  $\tilde{\iota}_1 \circ \tilde{\varphi}(R)$ is contained in the stratum  of $\tilde{\Acal}_1$ corresponding to $\relint(\Deltabar_1)$. This correspondence is described in \cite[Proposition 4.8]{gubler3}, showing also that $W_1^\circ:=\tilde{\iota}_1' \circ \tilde{\Phi}(R)$ is contained in the stratum $Z_{\Delta_1}$ of $\tilde{\Ecal}_1$ corresponding to $\relint(\Delta_1)$ for a suitable polytope $\Delta_1$ of $\rdop^n$ with image $\Deltabar_1$ in $\rtor$. By \cite[Remark 4.9]{gubler3}, this stratum is a torsor $\tilde{q}_1: Z_{\Delta_1} \rightarrow \tilde{\Bcal}$ with fibres isomorphic to a torus of dimension equal to $\codim(\Delta_1)$. Since $\tilde{\Phi}(R)$ is dense in $W$, it follows that $W_1^\circ$ is dense in $W_1$. We conclude that $W_1^\circ$ is contained in a fibre bundle over $\tilde{q}_1(W_1^\circ)$ with $\codim(\Delta_1)$-dimensional fibres. This and \eqref{equation 2} yield
\addtocounter{Claim}{1}
\begin{equation} \label{equation 6}
\dim S \geq \dim  \tilde{q}_0(\tilde{\Phi}_0(S)) = \dim \tilde{q}_1(W_1^\circ) \geq \dim W_1 - \codim \Delta_1.
\end{equation}
Since $W_1$ is an irreducible component of $\tilde{\Ycal}_1$, we have $\dim W_1 =d$. By \eqref{equation 5}, we get
$$ \dim W_1 - \codim \Delta_1 \geq d - \dim \Delta_S = \dim S.$$
We conclude that equality occurs everywhere in \eqref{equation 6} proving
\addtocounter{Claim}{1}
\begin{equation} \label{equation 7}
 \dim S = \dim  \tilde{q}_0(\tilde{\Phi}_0(S)).
\end{equation}
By \eqref{equation 4} and \eqref{equation 7}, the canonical simplex $\Delta_S$ is non-degenerate with respect to $f$. Using \eqref{equation 3} and \eqref{equation 7}, we conclude that $S$ is a $b(X)$-dimensional stratum of $\tilde{\Xcal}'$.
\end{Pf}

\begin{art} \label{proof of canonical measure} \rm
It remains to proof Theorem \ref{canonical measure}. We choose a strictly semistable alteration $\varphi_0:\Xcal' \rightarrow \Xcal_0$ as in \ref{setup} with generic fibre $f:(X')^{\rm an} \rightarrow \Xan$. Moreover, we may assume that the restriction of $f$ to $\Delta_S$ is a linear isomorphism onto a rational simplex of the canonical subset $S_X$ for all canonical simplices $\Delta_S$ of $S(\Xcal')$ which are non-degenerate with respect to $f$ (see the proof of \cite[Theorem 6.12]{gubler3}). We number these simplices of $T$ by $\Delta_1, \dots , \Delta_N$. 
By projection formula (\cite[Proposition 3.8]{gubler3}), we have 
$$f_*\left(c_1(f^*\overline{L_1}|_X) \wedge \dots \wedge c_1(f^*\overline{L_d}|_X)\right) = \deg(f) c_1(\overline{L_1}|_X) \wedge \dots \wedge c_1(\overline{L_d}|_X).$$
By \cite[Theorem 6.7 and Remark 6.8]{gubler3}, there are  numbers $r_S$ with
$$c_1(f^*\overline{L_1}|_X) \wedge \dots \wedge c_1(f^*\overline{L_d}|_X) = \sum_S r_S \delta_{\Delta_S}$$
where $S$ ranges over all strata of $\tilde{\Xcal}'$ such that the canonical simplex $\Delta_S$ of  the skeleton $S(\Xcal')$ is non-degenerate with respect to $f$. Note that the numbers $r_S$ are positive if all line bundles are ample.   This yields already properties (a)--(e) in Theorem \ref{canonical measure} and the last claim follows from Theorem \ref{semistable canonical measure}. \qed
\end{art}


\small{

}


\begin{thebibliography}
{99}
\bibitem{berkovich1}
V.G. Berkovich, 
Spectral theory and analytic geometry over nonarchimedean
fields. Mathematical Surveys and Monographs, 33. Providence,
RI: AMS (1990).
\bibitem{berkovich2}
V.G. Berkovich, \'Etale cohomology for non-archimedean 
analytic spaces.
Publ. Math. IHES 78 (1993), 5--161.
\bibitem{berkovich3}
V.G. Berkovich, 
Vanishing cycles for formal schemes. Invent. Math. 115-3 (1994), 
539--571.
\bibitem{berkovich4}
V.G. Berkovich, Smooth p-adic 
analytic spaces are locally contractible.
Invent. Math. 137-1 (1999), 1--84.
\bibitem{bosch}
S. Bosch,
Rigid analytische Gruppen mit guter Reduktion,
Math. Ann. 233 (1976), 193--205.
\bibitem{BL}
S. Bosch and W. L\"utkebohmert,
Degenerating abelian varieties,
Topology 30 (1991), 653--698.
\bibitem{chambert-loir}
A. Chambert-Loir, 
Mesure et \'equidistribution sur les espaces de Berkovich,
J. Reine Angew. Math. 595 (2006), 215--235.
\bibitem{cinkir}
Z. Cinkir, 
Zhang's Conjecture and the Effective Bogomolov Conjecture 
over function fields,
Invent. Math. 183 (2011), 517--562.
\bibitem{ega28}
A. Grothendieck,
\'El\'ements de g\'eom\'etrie alg\'ebrique IV
\'Etude locale des sch\'emas et des morphismes de sch\'emas III,
I.H.E.S.publ.Math. 28 (1966).
\bibitem{gubler1}
W. Gubler,
Tropical varieties for non-archimedean analytic spaces,
Invent. Math. 169 (2007), 321--376.
\bibitem{gubler2}
W. Gubler,
The Bogomolov conjecture for totally degenerate abelian varieties,
Invent. Math. 169 (2007), 377--400.
\bibitem{gubler4}
W. Gubler,
Equidistribution over function fields,
manuscripta math. 127 (2008), 485--510.
\bibitem{gubler3}
W. Gubler,
Non-archimedean canonical measures on abelian varieties,
Compositio Math. 146 (2010), 683--730.
\bibitem{dejong}
A. J. de Jong, Smoothness, 
semi-stability and alterations. Publ. Math.
IHES 83 (1996), 51--93.
\bibitem{lang1}
S. Lang,
Abelian varieties,
Springer-Verlag (1983).
\bibitem{lang2}
S. Lang,
Fundamentals of Diophantine Geometry,
Springer-Verlag (1983).
\bibitem{moriwaki3}
A. Moriwaki,
Bogomolov conjecture over function fields for stable curves 
with only irreducible fibers, 
Comp. Math., 105 (1997), 125--140.
\bibitem{moriwaki5}
A. Moriwaki,
Arithmetic height functions
over finitely generated fields,
Invent. Math. 140 (2000), 101--142.
\bibitem{ullmo}
E. Ullmo,
Positivit\'e et discr\'etion des points alg\'ebriques des courbes,
Ann. of Math. 147 (1998), 167--179.
\bibitem{yamaki1}
K. Yamaki,
Geometric Bogomolov's conjecture
for curves of genus $3$
over function fields,
J. Math. Kyoto Univ. 42 (2002),
57--81.
\bibitem{yamaki2}
K. Yamaki,
Effective calculation of the geometric height
and the Bogomolov conjecture
for hyperelliptic curves over
function fields,
J. Math. Kyoto. Univ. 48 (2008),
401--443.
\bibitem{zhang1}
S. Zhang,
Admissible pairing on a curve, 
Invent. Math. 112 (1993), 171--193.
\bibitem{zhang2}
S. Zhang,
Equidistribution of small points on abelian varieties,
Ann. of Math. 147 (1998), 159--165.
\end{thebibliography}
\end{document}